\documentclass[11pt]{article}
\usepackage[top=1in, bottom=1in, left=1in, right=1in]{geometry}
\usepackage{parskip}
\usepackage{amsmath, amssymb, amsfonts, amsthm, bm, mathrsfs}

\usepackage{graphicx}
\usepackage[percent]{overpic}
\usepackage{subcaption}

\usepackage{tikz}
\usetikzlibrary{positioning}

\tikzset{
  dim_node/.style={
    draw=blue!80,
    fill=blue!5, 
    rounded corners,
    align=center,
    inner sep=2pt
  },
  ht_tree/.style={
    grow=down,
    level distance=12mm,
    edge from parent/.style={draw, thick, ->},
  }
}

\usepackage{algorithm}
\usepackage{algpseudocode}

\usepackage{url}
\usepackage[colorlinks=true,
            linkcolor=blue,
            citecolor=blue,
            urlcolor=blue]{hyperref}
\usepackage[capitalize,nameinlink]{cleveref}

\usepackage{appendix}

\theoremstyle{plain}
\newtheorem{theorem}{Theorem}[section]
\newtheorem{lemma}[theorem]{Lemma}
\newtheorem{proposition}[theorem]{Proposition}

\theoremstyle{definition}
\newtheorem{definition}[theorem]{Definition}

\theoremstyle{remark}
\newtheorem{remark}[theorem]{Remark}

\title{Multilevel Adaptive-Rank Methods for Linear and Nonlinear Systems in the Hierarchical Tucker Format
}

\author{
William A. Sands\thanks{Theoretical Division, Los Alamos National Laboratory, Los Alamos, NM, 87545, USA (\href{mailto:wsands@lanl.gov}{\texttt{wsands@lanl.gov}})}
\and
Pierson T. Guthrey\thanks{Lawrence Livermore National Laboratory, Livermore, CA, 94550, USA (\href{mailto:guthrey1@llnl.gov}{\texttt{guthrey1@llnl.gov}})}
\and
Nathan V. Roberts\thanks{Sandia National Laboratories, Albuquerque, NM, 87123, USA (\href{mailto:nvrober@sandia.gov}{\texttt{nvrober@sandia.gov}})}
}

\date{}

\begin{document}

\maketitle

\begin{abstract}
We develop multilevel adaptive-rank iterative methods for the solution of linear and nonlinear systems arising from high-dimensional partial differential equations. Our contributions are threefold. First, we extend the projection method of Ballani and Grasedyck \cite{ballani2013projection} to enable flexible preconditioning of high-dimensional linear systems in low-rank tensor formats. Second, we construct multilevel preconditioning strategies by adapting geometric multigrid methods to the low-rank setting. In contrast to prior work, which primarily employs multigrid as a standalone solver, we emphasize its role as an efficient and robust preconditioner. Third, we integrate these techniques within an inexact Newton framework for the solution of nonlinear systems. The proposed methods are evaluated on a range of model problems, including both linear and nonlinear equations, to assess their convergence behavior and computational efficiency. The results demonstrate that multilevel adaptive-rank strategies yield robust and scalable preconditioners, providing effective solvers for high-dimensional problems in low-rank formats.
\end{abstract}


\section{Introduction}
\label{sec:introduction}

The solution of large-scale linear and nonlinear systems is a foundational problem in computational science and engineering that arises across a wide range of scientific applications. In this work, we are motivated by discretizations of partial differential equations (PDEs) that exhibit numerical stiffness, particularly in applications from kinetic theory. Representative examples include collisional models for plasmas and neutral gases, as well as phase-field models that describe interfacial dynamics in multiphase flows. Despite their differences, these systems are inherently multiscale and give rise to rich dynamics that are challenging to resolve in high-dimensional simulations.

A major obstacle in kinetic simulations is the intrinsic high dimensionality of phase space, which severely limits the applicability of conventional full-grid methods due to the curse of dimensionality. Recent advances in low-rank tensor methods (more broadly, tensor networks) have shown considerable promise in alleviating these challenges while preserving accuracy \cite{einkemmer2025review}. The key idea underlying these methods is to approximate the solution of a PDE in a separable, SVD-like format. When the number of terms in this representation is small, both storage requirements and computational complexity are significantly reduced. In time-dependent problems, this representation is evolved within the PDE, and algorithms such as the SVD or the high-order SVD (HOSVD) \cite{DeLathauwerHOSVD2000} are used to compress the solution by discarding components associated with small singular values. This paradigm has been successfully applied to a wide range of problems, including Vlasov-type plasma models \cite{kormann2015semi,allmann2022parallel,ye2022qtt_vp_system,ye2024qtt_vm_system,zheng2025HTACA}, transport equations \cite{sands2025transport,deshpande2026multigroupthermalradiationtransport}, and Fokker--Planck and diffusion equations \cite{dolgov2012qtt_fokker_planck,rodgers2023implicitHT,hamad2024krylov,nakao2026rail3d,hamad2025sylvester_preconditioned,wang2026dynamical_tt}.

Methods for linear systems in the low-rank setting have been developed along several research directions. Early work introduced Krylov subspace methods in tensor formats, such as the tensor-train (TT) GMRES (TT-GMRES) method \cite{dolgov2011tt_gmres}, which provides an iterative solution to high-dimensional linear systems in the TT format. Dolgov's approach is motivated by theory for inexact Krylov methods and uses a relaxed truncation tolerance that varies inversely with the residual norm. A limitation of this approach is its reliance on orthogonality. Since rounding cannot be applied after each tensor addition, the intermediate tensor ranks may grow significantly. These ideas were subsequently applied to time-dependent problems, including Fokker--Planck equations \cite{dolgov2012qtt_fokker_planck}, where low-rank structure is exploited throughout the time integration. Ballani and Grasedyck \cite{ballani2013projection} proposed an iterative method in the hierarchical Tucker (HT) format that projects the residual to low-dimensional subspaces, in a manner similar in spirit to Krylov methods. A distinguishing feature of their approach is that it does not require mutual orthogonality of the basis vectors, which is generally not preserved under truncation in low-rank formats. Coulaud et al.~\cite{coulaud2022robust} recently proposed a modification of Dolgov’s relaxed TT-GMRES algorithm that employs a backward error estimate as the stopping criterion. While this termination measure is simple to implement, it is not fundamentally more reliable than the relative residual, since it is based on the same inexact residual and relies on a linear perturbation interpretation that may break down in the presence of truncation.

Another class of iterative solvers is based on alternating minimization techniques. For example, the alternating minimal energy (AMEn) method \cite{dolgov2014amen} computes the solution of the global system by solving a sequence of smaller factor-wise problems. While these methods have been shown to be effective for certain structured problems, their convergence guarantees are limited in general settings and often depend on problem-specific properties. Although individual AMEn sweeps can be performed efficiently in high-dimensional problems, many sweeps may be required to achieve high accuracy, which can lead to significant computational cost. There has also been considerable interest in the development of direct methods based on algorithms for Sylvester-type matrix equations \cite{simoncini2020numerical}, which have been applied to three-dimensional convection--diffusion problems \cite{hamad2025sylvester_preconditioned,nakao2026rail3d}. While effective, these approaches often rely on problem-specific structure and recursive constructions that are difficult to generalize to high-dimensional settings.

Multigrid methods \cite{briggs2000multigrid,trottenberg2000multigrid} are a powerful class of multilevel iterative methods for linear systems that use a hierarchy of operators to target different frequency components of the error. When used as solvers or preconditioners, they are known to exhibit optimal or near-optimal convergence rates for elliptic and parabolic PDEs. However, their direct application to high-dimensional problems is limited by the curse of dimensionality. The tensor product structure inherent in low-rank representations greatly simplifies the construction of multigrid hierarchies and intergrid operators, making multilevel methods practical for high-dimensional problems. This idea was explored by Khoromskij and Khoromskaia \cite{khoromskij2009multigrid} to accelerate alternating least-squares methods, and later by Ballani and Grasedyck \cite{ballani2013projection}, who used multigrid to accelerate an iterative method. Andreev and Tobler \cite{andreev2015multilevel} proposed a space-time Petrov-Galerkin discretization of the linear diffusion equation in the HT format equipped with a multilevel preconditioner. Additional references can be found in the brief survey by Hackbusch \cite{hackbusch2015highDsystem}. Recent work by Grasedyck et al.~\cite{grasedyck2020parameter,grasedyck2023computing,grasedyck2025operator} extended these ideas to elliptic problems with variable coefficients and the design of effective smoothers.

In contrast to linear systems, the treatment of nonlinear problems in the low-rank setting remains relatively limited. Appelo and Cheng \cite{appelo2025lraa} employed Anderson acceleration to enhance fixed-point iterations by combining previous iterates to minimize the nonlinear residual. While straightforward to implement in low-rank formats, its convergence is generally weaker than Newton-type methods in strongly nonlinear regimes. Rodgers and Venturi \cite{rodgers2023implicitHT} extended Dolgov's TT-GMRES method \cite{dolgov2011tt_gmres} to the HT format and combined it with an inexact Newton method, resulting in one of the first approaches for nonlinear problems in this setting. Their method assumes that the action of the Jacobian can be obtained from an explicit linearization of the nonlinear problem. The linearization can be performed for problems with relatively simple nonlinearities, but it may be difficult or impractical for problems with strongly coupled nonlinear terms. Although convergence analysis is provided, their results do not offer insight into the computational cost of the iteration. Similarly, Adak et al.~\cite{adak2025tensor} developed inexact Newton methods for the TT format using a similar strategy, but provided limited information on the cost and scalability of the inner linear solves, particularly for problems with large mode sizes.

There remains a lack of robust and general-purpose preconditioning strategies for the linear systems arising in low-rank formulations, particularly in the context of nonlinear solvers. This work aims to address this gap by developing multilevel preconditioning strategies based on geometric multigrid methods in the low-rank setting. We combine these techniques with a flexible GMRES (FGMRES) method inspired by \cite{ballani2013projection} and inexact Newton methods to enable the efficient solution of both linear and nonlinear high-dimensional problems. In contrast to the Newton methods of Rodgers and Venturi \cite{rodgers2023implicitHT}, which assume an available linearization, we also consider approximate Jacobians obtained from finite differences of the residual in the spirit of Jacobian-free methods \cite{knoll2004jfnk}. This allows one to treat problems for which an explicit Jacobian is not available. When assessing the performance of the proposed methods, we provide details of the iteration data to demonstrate their robustness and scalability.

The remainder of this manuscript is organized as follows. In Section~\ref{sec:ht}, we provide a brief introduction to the HT format. While we are primarily interested in the HT format, the proposed methods can easily be extended to other tensor formats, such as Tucker \cite{tucker1966some,DeLathauwerHOSVD2000} and TT decompositions \cite{MPS2007,TToseledets2011}, provided that basic operations such as addition, scalar multiplication, truncation, inner products, and norm evaluation are available. Section~\ref{sec:linear solvers} presents the proposed methods for high-dimensional linear systems, including flexible GMRES and multigrid preconditioning strategies. These methods are extended to nonlinear problems using inexact Newton techniques in Section~\ref{sec:nonlinear solvers}. Numerical results are presented in Section~\ref{sec:numerics}, followed by conclusions and future work in Section~\ref{sec:conclusion}.

\section{HT Format: Notation and Preliminaries}
\label{sec:ht}


One of the primary challenges of high-dimensional systems is the exponential growth in storage requirements associated with standard storage containers, such as multidimensional arrays. To address the high dimensionality, we instead seek a representation for the components of the system in the HT format \cite{HackbuschKuhnHT2009}. The aim of this section is to introduce notation and provide a minimal introduction to the components of the format. For further details, we recommend consulting the work of Grasedyck \cite{grasedyck2010hierarchical}, which provides an excellent description of the format and its approximation properties, as well as implementation details of algorithms related to its construction.

The HT format constructs a hierarchical low-rank representation of a tensor by recursively partitioning its modes. This process is encoded in a binary structure called the \emph{dimension tree}, which defines the relationships between subsets of the modes. A formal definition for this can be stated as follows.

\begin{definition}[Dimension tree]
A binary tree $\mathcal{T}$, with each node represented by a subset of dimension labels $\{1,2,\ldots,d\}$, is called a \textit{dimension tree} if it satisfies the following conditions: 
\begin{enumerate}
    \item The root node $t_{\text{root}}$ contains all the dimensions $\{1,2,\ldots,d\}.$
    \item The two child nodes of every parent node are disjoint.
    \item Each leaf node contains a single index.
    \item Each parent node is the union of its two child nodes.
\end{enumerate}
\end{definition}

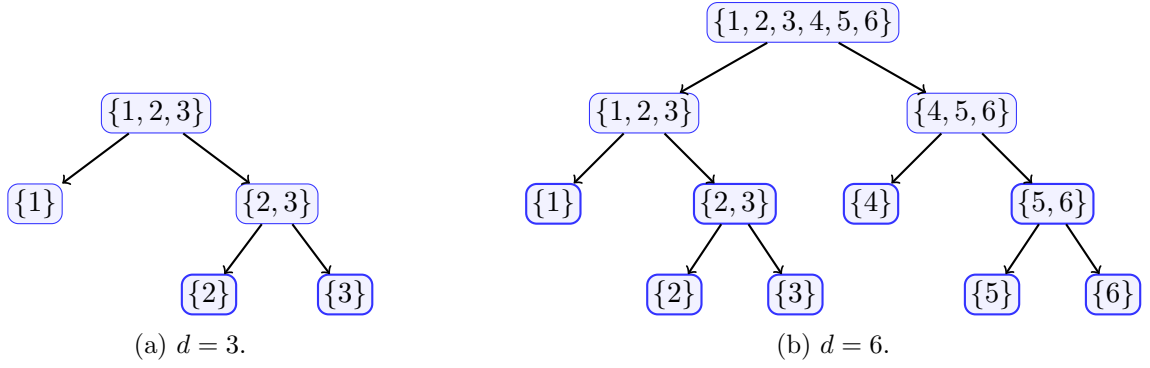
\begin{figure}[t]
\centering

\begin{subfigure}{0.45\linewidth}
\centering
\begin{tikzpicture}[ht_tree,
  level 1/.style={sibling distance=32mm},
  level 2/.style={sibling distance=18mm}]
  
\node[dim_node] {$\{1,2,3\}$}
  child { node[dim_node] {$\{1\}$} }
  child { node[dim_node] {$\{2,3\}$}
    child { node[dim_node] {$\{2\}$} }
    child { node[dim_node] {$\{3\}$} }
  };

\end{tikzpicture}
\caption{$d=3$.}
\end{subfigure}
\hfill
\begin{subfigure}{0.52\linewidth}
\centering
\begin{tikzpicture}[ht_tree,
  level 1/.style={sibling distance=42mm},
  level 2/.style={sibling distance=24mm},
  level 3/.style={sibling distance=16mm}]
  
\node[dim_node] {$\{1,2,3,4,5,6\}$}
  child { node[dim_node] {$\{1,2,3\}$}
    child { node[dim_node] {$\{1\}$} }
    child { node[dim_node] {$\{2,3\}$}
      child { node[dim_node] {$\{2\}$} }
      child { node[dim_node] {$\{3\}$} }
    }
  }
  child { node[dim_node] {$\{4,5,6\}$}
    child { node[dim_node] {$\{4\}$} }
    child { node[dim_node] {$\{5,6\}$}
      child { node[dim_node] {$\{5\}$} }
      child { node[dim_node] {$\{6\}$} }
    }
  };

\end{tikzpicture}
\caption{$d=6$.}
\end{subfigure}

\caption{Examples of dimension trees used in the HT format using $d=3$ (left) and $d=6$ (right).}
\label{fig:ht_dimension_trees}
\end{figure}

In what follows, we denote the set of leaf and non-leaf nodes as $\mathcal{L}(\mathcal{T})$ and $\mathcal{N}(\mathcal{T}) = \mathcal{T} \setminus \mathcal{L}(\mathcal{T})$, respectively. For each non-leaf node $\alpha \in \mathcal{N}(\mathcal{T})$, we denote its left and right children by $\alpha_l$ and $\alpha_r$, respectively, and the level of a node is its distance from the root node. We provide two examples of dimension trees in Figure \ref{fig:ht_dimension_trees} corresponding to $d = 3$ and $d = 6$.

A tensor $\mathcal{X} \in \mathbb{C}^{n_1 \times \cdots \times n_d}$ contains $d$ different modes, and we denote the mode size along mode $\mu$ by $n_{\mu}$. In addition, the entries along a mode $\mu$ are associated with an index set $\mathbb{I}_{\mu} = \{1,2, \cdots, n_{\mu}\}$. For brevity, we denote the tensor multi-index set as $\mathbb{I} = \mathbb{I}_{1} \times \mathbb{I}_{2} \times \cdots \times \mathbb{I}_{d}$. More generally, an index set can also be formed from any subset of the dimensions. In particular, for any subset $\alpha \subset \{1,\dots,d\}$ and its complement $\alpha^{c} = \{1,\dots,d\} \setminus \alpha$, we define
\begin{equation*}
    \mathbb{I}_\alpha := \prod_{\mu \in \alpha} \mathbb{I}_\mu, \quad \mathbb{I}_{\alpha^{c}} := \prod_{\mu \notin \alpha} \mathbb{I}_\mu.
\end{equation*}
With this notation, we define the process of matricization, or unfolding, along a set of dimensions as follows:
\begin{definition}[Matricization]
For $\mathcal{X} \in \mathbb{C}^{n_1 \times \cdots \times n_d}$ and 
$\alpha \subset \{1,\dots,d\}$, the mode-$\alpha$ matricization (or unfolding), denoted by
$X^{\alpha}$, is the matrix formed by merging the indices in 
$\alpha$ into a single row index $i^{\alpha} \in \mathbb{I}_{\alpha}$ and those in the complement $\alpha^c$ into a 
single column index $i^{\alpha^{c}} \in \mathbb{I}_{\alpha^{c}}$ such that
\begin{equation*}
    \mathcal{X}(i_{1}, i_{2}, \cdots, i_{d}) = X^{\alpha}\left(i^{\alpha}, i^{\alpha^{c}}\right).
\end{equation*}
\end{definition}
\noindent Note that \emph{vectorization} is a special case in which one performs matricization along \emph{all} dimensions of the tensor. Similarly, one can define a \emph{dematricization} (or refolding) of a matrix into a tensor with compatible dimensions by reversing the matricization process. It is important to note that the matricization (and similarly dematricization) depends on the ordering of the modes of the tensor, which we assume to be sorted in ascending order.

The dimension tree in the HT format defines the hierarchical grouping of modes and their associated topology. The latter is critical to the construction of the basis and the so-called transfer tensors in the decomposition. In particular, given a tensor $\mathcal{X}$, we have the following recursive property. For each non-leaf node $\alpha \in \mathcal{N}(\mathcal{T})$ with children $\alpha_l$ and $\alpha_r$, the column space of $X^{\alpha}$ lies within that of the Kronecker product $X^{\alpha_r} \otimes X^{\alpha_l}$. This property is formally stated below.

\begin{lemma}[Grasedyck {\cite{grasedyck2010hierarchical}}]
Let $\mathcal{X} \in \mathbb{C}^{n_1 \times \cdots \times n_d}$ 
and let $\mathcal{T}$ be a dimension tree.
For any non-leaf node $\alpha = \alpha_l \cup \alpha_r \in \mathcal{T}$,
the column space of the mode-$\alpha$ matricization satisfies
\[
\mathrm{range}\!\left( X^{\alpha} \right)
\subseteq
\mathrm{range}\!\left(
X^{\alpha_r} \otimes X^{\alpha_l}
\right).
\]
\label{lem:HTD_foundation}
\end{lemma}

For tensors in the HT format, there is a concept of hierarchical rank which is defined as follows.
\begin{definition}[Hierarchical rank]
Given a dimension tree \(\mathcal{T}\), the \textit{hierarchical rank} of a tensor \(\mathcal{X}\) is defined as
\[
(r_{\alpha})_{\alpha\in\mathcal{T}} = (\mathrm{rank}(X^{\alpha}))_{\alpha\in\mathcal{T}}.
\]
We denote the set of HT tensors with hierarchical ranks bounded by \((r_\alpha)_{\alpha\in\mathcal{T}}\) by
\[
\mathcal{H}\text{-Tucker}((r_\alpha)_{\alpha\in\mathcal{T}}) = \left\{ \mathcal{X} \in \mathbb{C}^{n_1\times \cdots \times n_d} \,\middle|\, \mathrm{rank}(X^{\alpha}) \leq r_\alpha,\ \forall \alpha \in \mathcal{T} \right\}.
\]
\end{definition}
\noindent The hierarchical rank \((r_\alpha)_{\alpha\in\mathcal{T}}\) consists of \(2d-1\) values, one for each node of the binary tree, and the rank at the root node of the tree is fixed to 1 by convention.

\begin{definition}[Transfer tensor]
Let $\mathcal{X}$ be a $d$-dimensional tensor and let $\mathcal{T}$ be a dimension tree.
Consider a non-leaf node $\alpha \in \mathcal{T}$ with children
$\alpha_l$ and $\alpha_r$. Let 
\[
U_{\alpha} \in \mathbb{C}^{\mathbb{I}_{\alpha} \times r_{\alpha}}, \quad
U_{\alpha_l} \in \mathbb{C}^{\mathbb{I}_{\alpha_l} \times r_{\alpha_l}}, \quad
U_{\alpha_r} \in \mathbb{C}^{\mathbb{I}_{\alpha_r} \times r_{\alpha_r}}
\]
be basis matrices whose columns span the ranges of
$X^{\alpha}$, $X^{\alpha_l}$, and $X^{\alpha_r}$, respectively. Then, by Lemma \ref{lem:HTD_foundation} there exists a matrix
$$
B_{\alpha} \in
\mathbb{C}^{(r_{\alpha_l} r_{\alpha_r}) \times r_{\alpha}},
$$
such that
\begin{equation*}
    U_{\alpha}
    =
    \left(
    U_{\alpha_r}
    \otimes
    U_{\alpha_l}
    \right)
    B_{\alpha}.
\end{equation*}
Reshaping $B_{\alpha}$ into a third-order tensor
$$
\mathcal{B}_{\alpha}
\in
\mathbb{C}^{r_{\alpha_l} \times r_{\alpha_r} \times r_{\alpha}}
$$
yields the \emph{transfer tensor} associated with node $\alpha$.
\label{def:transfer tensors}
\end{definition}

By unifying these concepts, we define the HT decomposition (HTD) of a tensor as follows.
\begin{definition}[HTD]
\label{def:HTD}  
Given \(\mathcal{X} \in \mathcal{H}\text{-Tucker}((r_\alpha)_{\alpha \in \mathcal{T}})\), the HTD of \(\mathcal{X}\) consists of a set of transfer tensors \((\mathcal{B}_\alpha)_{\alpha \in \mathcal{N}(\mathcal{T})}\) and leaf node matrices \((U_\alpha)_{\alpha \in \mathcal{L}(\mathcal{T})}\), such that the recursive property in Definition \ref{def:transfer tensors} holds for all non-leaf nodes.  
We denote this decomposition as  
$$
\mathcal{X} = ((\mathcal{B}_\alpha)_{\alpha\in\mathcal{N}(\mathcal{T})},\ (U_\alpha)_{\alpha\in\mathcal{L}(\mathcal{T})}).
$$
\end{definition}

As seen in Definition \ref{def:HTD}, the HTD requires storing the transfer tensors and basis matrices, respectively, at the non-leaf and leaf nodes of the dimension tree. The following proposition relates the dimension, mode sizes, and hierarchical rank with the storage complexity of this format.
\begin{proposition}\label{prop:HTD_storage}
Let $n = \max\limits_{1 \leq \mu \leq d} \{n_{\mu}\}$ and $r = \max\limits_{\alpha \in \mathcal{T}} \{r_\alpha\}$. Then the storage complexity of a tensor $\mathcal{X} \in \mathbb{C}^{n_1 \times \cdots \times n_d}$ in the HT format with hierarchical rank $(r_{\alpha})_{\alpha\in\mathcal{T}}$ is at most
\begin{equation*}
    dnr + (d - 2)r^3 + r^2.
\end{equation*}
\end{proposition} 
\noindent Note that the dimension $d$ no longer appears in the exponent of the mode size $n$. Consequently, the condition $r \ll n$ is sufficient but not necessary for this format to avoid the curse of dimensionality.

Later on we shall measure the effectiveness of the proposed adaptive-rank methods in compressing the solutions to systems. For this, we define the compression rate as follows:
\begin{definition}[Compression rate]
\label{def:compression_rate}  
Let $\mathcal{X} \in \mathbb{C}^{n_1 \times \cdots \times n_d}$ and let
$\widehat{\mathcal{X}}$ denote its HT representation with (adaptive)
hierarchical ranks $(r_\alpha)_{\alpha\in\mathcal{T}}$.
We define the \emph{compression rate} as the ratio of the degrees of
freedom in the full tensor representation to the degrees of freedom in
the HT representation:
$$
\mathrm{Compression~Rate~}(\widehat{\mathcal{X}}) :=
\frac{\texttt{ndofs}(\mathcal{X})}{\texttt{ndofs}(\widehat{\mathcal{X}})}.
$$
\end{definition}

\begin{remark}[Asymptotic compression]
\label{rem:asymptotic_compression}  
Consider the case where $n_{\mu} \equiv n$, for $\mu = 1, 2, \cdots,d$, and assume that the hierarchical ranks are bounded such that $r_\alpha \lesssim r$ with $r \ll n$. Then, the compression rate in $d$ dimensions for the HT format scales roughly as
$$
\mathrm{Compression~Rate~} (\widehat{\mathcal{X}}) = \mathcal{O} \left(n^{d-1}\right).
$$
\end{remark}
\noindent This simple calculation highlights the potential of hierarchical low-rank formats in dramatically reducing the storage costs associated with more conventional storage formats.

When performing operations such as basic arithmetic on tensors, the use of truncation or rounding is crucial to maintaining tractability. This truncation step introduces a controlled approximation error with each application. In the HT format, this approximation is obtained by truncating singular components in the matricizations associated with the nodes of the dimension tree, following the HOSVD framework \cite{DeLathauwerHOSVD2000}. The resulting error can be bounded in terms of the discarded singular-value tails accumulated over the tree. This yields the following approximation property.
\begin{theorem}[\cite{grasedyck2010hierarchical}, Theorem 3.11]
    \label{thm:HTD truncation}
    Let $\mathcal{X} \in \mathbb{C}^{n_{1} \times n_{2} \times \cdots \times n_{d}}$ and let $\widehat{\mathcal{X}}$ denote its truncation to an HTD with dimension tree $\mathcal{T}$ and hierarchical ranks $\left\{ r_{t} : t\in \mathcal{T} \right\}$. Then the resulting approximation satisfies
    \begin{equation*}
        \left\lvert\left\lvert \mathcal{X} - \widehat{\mathcal{X}} \right\rvert\right\rvert_{F} \leq \sqrt{\sum_{t\in \mathcal{T}}\sum_{i > r_{t}} \sigma_{t,i}^{2} } \leq \sqrt{2d-3} \left\lvert\left\lvert \mathcal{X} - \mathcal{X}_{\text{opt}} \right\rvert\right\rvert_{F},
    \end{equation*}
    where $\mathcal{X}_{\text{opt}}$ denotes the best approximation of $\mathcal{X}$ in $\mathcal{H}\text{-Tucker}((r_\alpha)_{\alpha \in \mathcal{T}})$. Here $\sigma_{t,i} = \sigma_{i}\left( X^t \right)$ denotes the singular values of the mode-$t$ matricization of $\mathcal{X}$ and $ \left\lvert\left\lvert \cdot \right\rvert\right\rvert_{F} $ denotes the tensor Frobenius norm, which is defined as
    \begin{equation}
        \label{eq:frobenius norm}
        \left\lvert\left\lvert \mathcal{X} \right\rvert\right\rvert_{F} = \left( \sum_{i_{1} = 1}^{n_{1}}\cdots \sum_{i_{d} = 1}^{n_{d}} |\mathcal{X}_{i_{1,}\cdots, i_{d}}|^{2} \right)^{1/2}.
    \end{equation}
    Note that this is induced by the inner product
    \begin{equation}
        \label{eq:inner product}
        \left\langle \mathcal{X}, \mathcal{Y}\right\rangle = \sum_{i_{1} = 1}^{n_{1}}\cdots \sum_{i_{d} = 1}^{n_{d}} \bar{\mathcal{X}}_{i_{1,}\cdots, i_{d}} \mathcal{Y}_{i_{1,}\cdots, i_{d}},
    \end{equation}
    where the bar denotes elementwise conjugation.
\end{theorem}

\begin{remark}[\cite{kressner2014algorithm}, Remark 6.1] The hierarchical ranks $\left\{ r_{t} : t\in \mathcal{T} \right\}$ can be selected to meet prescribed truncation tolerances. In particular, given absolute and relative truncation tolerances $\epsilon_{\text{abs}}$ and $\epsilon_{\text{rel}}$, for any node $t\in \mathcal{T} \setminus \{t_{\text{root}}\}$ one can choose the rank $r_{t}$ so that
\begin{align}
    \left\lvert\left\lvert \mathcal{X} - \widehat{\mathcal{X}} \right\rvert\right\rvert_{F} &\leq \epsilon_{\text{abs}}, \quad \sqrt{\sum_{i > r_{t}} \sigma_{t,i}^{2}} \leq \frac{\epsilon_{\text{abs}}}{\sqrt{2d-3}}, \label{eq:abs trunc tolerance} \\
    \left\lvert\left\lvert \mathcal{X} - \widehat{\mathcal{X}} \right\rvert\right\rvert_{F} &\leq \epsilon_{\text{rel}} \left\lvert\left\lvert \mathcal{X} \right\rvert\right\rvert_{F}, \quad \sqrt{\sum_{i > r_{t}} \sigma_{t,i}^{2}} \leq \frac{\epsilon_{\text{rel}} \left\lvert\left\lvert \mathcal{X} \right\rvert\right\rvert_{F}}{\sqrt{2d-3}}, \label{eq:rel trunc tolerance}
 \end{align}
 simultaneously hold.
\end{remark}

We denote the process of truncation using the operator 
\begin{equation*}
    \widehat{\mathcal{X}} := \mathfrak{T}\left(\mathcal{X}\right),
\end{equation*}
which acts on a single operand; however, we shall generally omit the ``hat" on the output for convenience. Additionally, the methods proposed in this work frequently make use of truncated sums of tensors, which we denote by $\mathfrak{T}_{+}\left( \left\{ \mathcal{X}_{1}, \mathcal{X}_{2}, \cdots, \mathcal{X}_{s}\right\} \right)$.

The algorithms developed in this work rely on a small set of fundamental operations in the HT format, including orthogonalization, inner products, and truncation. Complexity estimates for these operations, along with implementation details, can be found in Kressner and Tobler \cite{kressner2014algorithm}; we briefly summarize them here for completeness. To simplify the analysis, these estimates assume that the HTD has a uniform mode size $n$, fixed dimension $d$, and maximum hierarchical rank $r$. Under these assumptions, orthogonalization, inner product and norm evaluation, and truncation of a single HTD each have cost $\mathcal{O}(dnr^{2} + dr^{4})$. For sums of HTDs, it is advantageous to interleave addition with truncation rather than forming the full sum prior to compression, as this avoids excessive intermediate storage growth. In particular, for a sum of $s$ HTDs with maximum hierarchical rank $r$ and mode size $n$, the total cost of the combined addition and truncation is $\mathcal{O}(dns^{2}r^{2} + ds^{2}r^{4} + ds^{3}r^{3})$. Finally, we note that inner products and norms should be computed using orthogonalized HTDs to ensure numerical stability and accuracy.

\section{Iterative Methods for Linear Systems}
\label{sec:linear solvers}

In this section we provide details regarding the iterative methods used to solve high-dimensional linear systems of equations. We begin by introducing the flexible GMRES (FGMRES) solver used in this work, which is based on the adaptive-rank projection method of Ballani and Grasedyck \cite{ballani2013projection}, and supports both left- and right-preconditioning. Once we have introduced the linear solver, we then describe the proposed multigrid preconditioning strategy, which is compatible with the low-rank format.

\subsection{The FGMRES Method for Linear Systems}
\label{subsec:fgmres}


Consider the high-dimensional linear system
\begin{equation}
    \label{eq:high-D linear system}
    \mathcal{A}x = b,
\end{equation}
where $\mathcal{A} \in \mathbb{C}^{\mathbb{I \times \mathbb{I}}}$ is a linear operator, $b \in \mathbb{C}^{\mathbb{I}}$ is the right-hand side, and $x \in \mathbb{C}^{\mathbb{I}}$ is the solution. In the applications considered in this work, the operator $\mathcal{A}$ is assumed to be available in a tensor product format given by 
\begin{equation}
    \label{eq:A as tensor product}
    \mathcal{A} = \sum_{k=1}^{r_A} \bigotimes_{\mu=1}^{d} A_{\mu,k}.
\end{equation}
A special case of this is the Kronecker sum
\begin{equation*}
    \mathcal{A} = A_{1} \oplus A_{2} \oplus \cdots \oplus A_{d} \equiv \sum_{\mu = 1}^{d} I_{1} \otimes \cdots \otimes I_{\mu-1} \otimes A_{\mu} \otimes I_{\mu+1} \otimes \cdots \otimes I_{d},
\end{equation*}
where $I_{\mu} \in \mathbb{C}^{n_{\mu} \times n_{\mu}}$ denotes the identity matrix. In practice, the factors of $\mathcal{A}$ arise from finite-difference or finite-element discretizations of PDEs. Notably, the methods developed here do not require explicit construction of these matrices. The essential requirement is that the user implement the action of $\mathcal{A}$ on a given tensor $x$. 

We seek an iterative solution of \eqref{eq:high-D linear system} using a projection-based method inspired by Krylov subspace techniques. The proposed method builds on the projection framework of Ballani and Grasedyck \cite{ballani2013projection} and provides modifications to support flexible preconditioning. A central component of this work is the development of multilevel preconditioners that preserve the tensor product structure of the system. While our primary focus is on right preconditioning, it is convenient to consider a two-sided preconditioned formulation. To this end, we introduce an auxiliary variable $u \in \mathbb{C}^{\mathbb{I}}$ such that
\begin{equation*}
    x = \mathcal{M}_{r}^{-1}u
\end{equation*}
and consider the equivalent system
\begin{equation*}
    \mathcal{M}_{l} \mathcal{A} \mathcal{M}_{r}^{-1}u = \mathcal{M}_{l}b.
\end{equation*}
Here, $\mathcal{M}_{l}$ and $\mathcal{M}_{r}^{-1}$ denote (possibly nonlinear and iteration-dependent) preconditioning operators. Although these operators are not required to be linear across iterations, we assume that each application admits a tensor product structure, which enables efficient dimension-wise application in low-rank formats. The residual associated with the two-sided preconditioned system is given by
\begin{equation*}
    r = \mathcal{M}_{l} \left( b - \mathcal{A} x \right),
\end{equation*}
and convergence is measured in the Frobenius norm \eqref{eq:frobenius norm}.

The method consists of outer and inner iteration phases. The outer iteration monitors the residual, while the inner iteration constructs a low-dimensional subspace over which the residual is minimized. To describe a single outer iteration, let $x_0 \in \mathbb{C}^{\mathbb{I}}$ be the current iterate and define
\begin{equation*}
    r_{0} = \mathcal{M}_{l} \left( b - \mathcal{A} x_{0} \right).
\end{equation*}
Next, we write $x = x_{0} + e$, where $e \in \mathbb{C}^{\mathbb{I}}$ represents the error. Introducing $\delta \in \mathbb{C}^{\mathbb{I}}$ such that
\begin{equation*}
    e = \mathcal{M}_r^{-1} \delta,
\end{equation*}
we obtain the correction equation
\begin{equation}
    \label{eq:correction equation}
    \mathcal{M}_{l} \mathcal{A} \mathcal{M}_{r}^{-1}\delta = r_{0}.
\end{equation}

The inner iteration seeks an approximate solution to \eqref{eq:correction equation} in a subspace of dimension $m$. This subspace is constructed via an Arnoldi-like process, initialized with $v_1 = r_0 / \|r_0\|$. At iteration $m$, we define
\begin{equation}
    \label{eq:test and trial vectors}
    z_j = \mathcal{M}_{r,j}^{-1} v_j, 
    \qquad
    w_j = \mathcal{M}_{l,j} \mathcal{A} z_j, 
    \qquad j = 1, \dots, m,
\end{equation}
where the preconditioners may vary across iterations, as indicated by the subscripts. Consequently, we store three sets of basis vectors during the inner iteration, namely
\[
V_m = [v_1, \dots, v_m] \in \mathbb{C}^{\mathbb{I} \times m}, \quad
Z_m = [z_1, \dots, z_m] \in \mathbb{C}^{\mathbb{I} \times m}, \quad
W_m = [w_1, \dots, w_m] \in \mathbb{C}^{\mathbb{I} \times m}.
\]
These basis vectors correspond to the residual, trial, and test spaces, respectively. 

At iteration $m$, we seek an approximation to the error $e$ of the form
\begin{equation*}
    e \approx Z_m y, \quad y \in \mathbb{C}^m,
\end{equation*}
which corresponds to approximating the residual as
\begin{equation*}
    r_{0} \approx W_{m}y.
\end{equation*}
Enforcing the Petrov--Galerkin condition
\begin{equation*}
    \left\langle W_{m}, r_{0} - W_{m} y \right\rangle = 0,
\end{equation*}
yields the $m \times m$ Gramian system
\begin{equation}
    \label{eq:y_new projection coefficients}
    W_{m}^{*} W_{m} y = W_{m}^{*} r_{0}.
\end{equation}
After solving for $y \in \mathbb{C}^m$, we update the residual and the solution according to
\begin{equation}
    \label{eq:r_new and x_new projection}
    r = r_0 - W_m y, \quad
    x = x_0 + Z_m y.
\end{equation}

Following Ballani and Grasedyck \cite{ballani2013projection}, the update is accepted only if a sufficient decrease in the residual is achieved:
\begin{equation}
    \label{eq:FGMRES inner acceptance criterion}
    \|r\|_{F} \leq \eta_{\text{rel}} \|r_0\|_{F}, \quad \eta_{\text{rel}} \in (0,1).
\end{equation}
If this condition is not satisfied, the subspace is enriched by increasing $m$. In order to obtain the next basis vector, we project $w_m$ onto the complement of $\mathrm{span}(V_m)$. Specifically, we compute $\alpha \in \mathbb{C}^m$ from the $m \times m$ Gramian system
\begin{equation}
    \label{eq:v_new projection alpha coefficients}
    V_m^* V_m \alpha = V_m^* w_m,
\end{equation}
and define
\begin{equation}
    \label{eq:v_new projection}
    v_{m+1} = w_m - V_m \alpha,
\end{equation}
which is followed by a normalization step. Upon termination of the inner iteration, the outer residual is updated and the process is restarted.

\subsection{Nonorthogonal FGMRES Method for Tensor Decompositions}
\label{subsec:fgmres tensor}


Although relatively simple, the method described in the previous section becomes impractical for high-dimensional problems due to the exponential growth of degrees of freedom. To address this issue, we incorporate low-rank tensor decompositions, which represent the solution and intermediate quantities in a compressed format.

In exact arithmetic, the method of the previous section corresponds to a special case of the FGMRES method proposed by Saad \cite{saad1993flexible}. However, coupling with low-rank tensor formats introduces several important modifications. In particular, basic tensor operations such as addition and subtraction lead to an increase in rank, and therefore require truncation to maintain computational tractability. 

The introduction of truncation has several important consequences. First, tensor truncation defines a nonlinear mapping, and therefore compositions of linear operations with truncation result in a nonlinear algorithm. Consequently, the Arnoldi process used to construct the basis vectors is no longer linear, and the resulting vectors are not mutually orthogonal. Therefore, the method no longer generates a true Krylov subspace of the form
\begin{equation*}
    \mathcal{K}_m(\mathcal{A}, r_0) = \mathrm{span}\left\{r_0, \mathcal{A}r_0, \mathcal{A}^2 r_0, \dots \mathcal{A}^{m-1} r_0 \right\}.
\end{equation*}

Instead, the proposed method constructs an approximate subspace through successive applications of operators that include effects from truncation. For example, the construction of new basis vectors for the residual, defined by \eqref{eq:v_new projection}, is modified to
\begin{equation*}
    v_{m+1} = \mathfrak{T}_{+}\left(\{w_m, -\alpha_{1} v_{1}, \cdots, -\alpha_{m}v_{m}\}\right).
\end{equation*}
The coefficients $\alpha$ are obtained from the Gramian system \eqref{eq:v_new projection alpha coefficients}, whose components are calculated through a series of inner products, namely
\begin{equation*}
    \left( V_{m}^{*}V_{m} \right)_{i,j} 
    = \left\langle v_{i}, v_{j} \right\rangle, 
    \quad 
    \left( V_{m}^{*}w_{m} \right)_{i} 
    = \left\langle v_{i}, w_{m} \right\rangle, 
    \quad i,j = 1,\dots,m.
\end{equation*}
Here the inner product $\langle\cdot,\cdot\rangle$ is defined according to \eqref{eq:inner product} and can be efficiently evaluated in the low-rank format. Due to truncation, which is interleaved with addition in the definition of $v_{m+1}$, this step only approximately removes components in $\mathrm{span}(V_m)$. Consequently, the vectors $\{v_j\}$ are not mutually orthogonal in general. This loss of orthogonality propagates to the construction of the Gramian system \eqref{eq:y_new projection coefficients}, as well as to the updates of the residual and solution defined by \eqref{eq:r_new and x_new projection}. In particular, the initial residual is defined by
\begin{equation*}
    r_{0} = \mathcal{M}_{l} \Big( \mathfrak{T} \left( b - \mathcal{A} x_{0} \right) \Big),
\end{equation*}
and is updated according to
\begin{equation*}
    r = \mathfrak{T}_{+} \left(\{r_0, - y_{1}w_{1}, \cdots, -y_{m}w_{m}\}\right), \quad
    x = \mathfrak{T}_{+} \left(\{x_0,   y_{1}z_{1}, \cdots,  y_{m}z_{m}\}\right).
\end{equation*}

Although not indicated explicitly, we assume that the operator $\mathcal{A}$ and the preconditioners $\mathcal{M}_{l,j}$ and $\mathcal{M}_{r,j}$ internally incorporate truncation to maintain a prescribed low-rank representation. As a result, the trial and test vectors defined in \eqref{eq:test and trial vectors} are also generated through controlled nonlinear mappings. These observations motivate the use of FGMRES, which does not require a fixed linear operator or a true Krylov structure. Instead, the method can be interpreted as a projection onto adaptively generated subspaces
\begin{equation*}
V_m = \mathrm{span}\{v_1, \dots, v_m\}, \quad Z_m = \mathrm{span}\{z_1, \dots, z_m\},
\end{equation*}
with residual minimization enforced in the corresponding test space
\begin{equation*}
W_m = \mathrm{span}\{w_1, \dots, w_m\}.
\end{equation*}

\subsection{Multigrid Methods}
\label{subsec:multigrid_ht}


Multigrid methods are a widely used and effective approach for solving linear systems with elliptic and parabolic structure. In the context of low-rank methods, they have primarily been used as standalone solvers for high-dimensional systems of the form \eqref{eq:high-D linear system} (see, e.g., \cite{ballani2013projection,grasedyck2020parameter,grasedyck2023computing,grasedyck2025operator}). In contrast, this work considers geometric multigrid (GMG) as a preconditioner for the iterative solver introduced in the previous subsections. This strategy exploits the tensor product structure of the problem to construct multilevel hierarchies that are compatible with iterative solvers in the low-rank format.

Let $\Omega = \Omega_{1} \times \Omega_{2} \times \cdots \times \Omega_{d}$ denote a $d$-dimensional tensor product domain. We consider a sequence of nested grids $\{\Omega^{(\ell)}\}_{\ell=1}^{L}$, where $\ell = 1$ corresponds to the coarsest level and $\ell = L$ to the finest level, i.e., $\Omega$. Let $h_{\mu}^{(\ell)}$ denote the uniform grid spacing at level $\ell$ in dimension $\mu$. We assume a standard coarsening strategy in which adjacent levels differ by at most a factor of two: 
\begin{equation*}
    h_{\mu}^{(\ell)} / h_{\mu}^{(\ell-1)} \leq 2, \quad \mu = 1,2, \cdots, d.    
\end{equation*}

The mesh hierarchy is constructed from a user-specified fine grid with a prescribed upper bound $n_{\mu,\min}$ on the mode sizes at the coarsest level. Let $n_\mu^{(\ell)}$ denote the mode size in dimension $\mu$ at level $\ell$. For the nested vertex-centered grids considered in this paper, we require $n_\mu^{(\ell)} = 2^{m} + 1$. We note that in the case of a cell-centered discretization, the analogous requirement is $n_\mu^{(\ell)} = 2^{m}$. For each dimension $\mu$, we define the number of active coarsening levels as
\begin{equation*}
    L_{\mu} := \min \left\{ \ell \in \mathbb{Z}_{+} : n_\mu^{(\ell)} \leq n_{\mu,\min} \right\},
\end{equation*}
and set the total number of levels to be
\begin{equation*}
L := \max_{\mu=1,\ldots,d} L_\mu.
\end{equation*}
This convention ensures that all dimensions share a common hierarchy depth, even when their mode sizes differ. While this assumption is not strictly necessary, the uniform construction it provides simplifies the bookkeeping associated with the hierarchy and the application of operators across its levels.

On each level $\ell$, we define a discretized operator of the form \eqref{eq:A as tensor product}:
\begin{equation}
    \label{eq:Level-wise operator A}
    \mathcal{A}^{(\ell)} 
    = \sum_{k=1}^{r_A} \bigotimes_{\mu=1}^{d} A_{\mu,k}^{(\ell)},
\end{equation}
where each $A_{\mu,k}^{(\ell)}$ is a one-dimensional operator that acts only along mode $\mu$, and $r_A$ denotes the rank of the operator.

To transfer information between levels, we introduce restriction and prolongation operators
\begin{align*}
\mathcal{R}^{(\ell)} &: \mathbb{C}^{\mathbb{I}^{(\ell)}} \to \mathbb{C}^{\mathbb{I}^{(\ell-1)}}, 
\quad \ell = 2,\cdots, L, \\
\mathcal{P}^{(\ell)} &: \mathbb{C}^{\mathbb{I}^{(\ell)}} \to \mathbb{C}^{\mathbb{I}^{(\ell+1)}},
\quad \ell = 1,2,\cdots, L-1.
\end{align*}
The restriction operator maps the grid data at a particular level $\ell$ to the next coarser level $\ell -1$, while prolongation interpolates from level $\ell$ to the next finer level $\ell + 1$. Owing to the tensor product structure of the mesh hierarchy, these operators admit separable representations
\begin{equation*}
    \mathcal{R}^{(\ell)} = \bigotimes_{\mu=1}^{d} R_{\mu}^{(\ell)}, \quad \mathcal{P}^{(\ell)} = \bigotimes_{\mu=1}^{d} P_{\mu}^{(\ell)}.
\end{equation*}
As a result, both restriction and prolongation can be applied in a dimension-wise manner, and, since they are rank-1 operators, they do not change the rank of the operand.

To account for anisotropic mode sizes, transfer schedules are constructed dimension by dimension starting from the finest level in the hierarchy. For each dimension $\mu$, the mode size is successively coarsened until a prescribed threshold is reached. To formalize this process, we introduce a Boolean restriction schedule, represented by a matrix with entries
\begin{equation*}
    \delta_\mu^{(\ell)} \in \{0,1\}, \quad \mu = 1,2, \cdots,d, \quad \ell = 2, 3, \cdots, L.
\end{equation*}
Here, $\delta_\mu^{(\ell)} = 1$ indicates that dimension $\mu$ is coarsened when transferring from level $\ell$ to level $\ell-1$, while $\delta_\mu^{(\ell)} = 0$ indicates that no coarsening is applied in dimension $\mu$ at that step. The mode-wise restriction operators are then defined as
\begin{equation*}
    R_\mu^{(\ell)} =
    \begin{cases}
    \widetilde{R}_{\mu}^{(\ell)}, & \delta_{\mu}^{(\ell)} = 1, \\
    I_{\mu}, & \delta_{\mu}^{(\ell)} = 0,
    \end{cases}
    \quad
    \ell = 2,3, \cdots, L
\end{equation*}
where $\widetilde{R}_{\mu}^{(\ell)}$ denotes the $n_{\mu}^{(\ell-1)} \times n_{\mu}^{(\ell)}$ restriction operator defined by the standard full-weighting stencil
\begin{equation*}
    u_i^{(\ell-1)} = \frac{1}{4}\,u_{2i-1}^{(\ell)} +  \frac{1}{2}\,u_{2i}^{(\ell)} +  \frac{1}{4}\,u_{2i+1}^{(\ell)},
\end{equation*}
with appropriate modifications to account for the boundary conditions. 

The prolongation schedule is obtained by shifting the restriction schedule by one level, so that refinement occurs exactly in dimensions that were coarsened on the next finer level:
\begin{equation*}
    P_\mu^{(\ell)} =
    \begin{cases}
    \widetilde{P}_{\mu}^{(\ell)}, & \delta_{\mu}^{(\ell+1)} = 1, \\
    I_{\mu}, & \delta_{\mu}^{(\ell+1)} = 0,
    \end{cases}
    \quad
    \ell = 1,2, \cdots, L-1
\end{equation*}
where $\widetilde{P}_{\mu}^{(\ell)}$ denotes the $n_{\mu}^{(\ell+1)} \times n_{\mu}^{(\ell)}$  prolongation operator. We use direct injection at shared grid points and linear interpolation at the off-grid points, which yields the respective stencils
\begin{equation*}
    u_{2i}^{(\ell+1)} = \,u_{i}^{(\ell)}, \quad u_{2i+1}^{(\ell+1)} = \frac{1}{2}\,u_{i}^{(\ell)} +  \frac{1}{2}\,u_{i+1}^{(\ell)}.
\end{equation*}
As with the restriction operator, these stencils are modified near the domain boundaries as required.

To reduce high-frequency error components during inter-level transfers, we apply an adaptive-rank Jacobi smoother, which decomposes the operator \eqref{eq:Level-wise operator A} as
\begin{equation*}
    \mathcal{A}^{(\ell)} = \mathcal{D}^{(\ell)} + \mathcal{N}^{(\ell)}, \quad \ell = 1, 2, \cdots, L,
\end{equation*}
where
\begin{equation*}
    \mathcal{D}^{(\ell)} := \sum_{k=1}^{r_{A}}\bigotimes_{\mu=1}^{d} \mathrm{diag}\left( A_{\mu, k}^{(\ell)} \right), \quad \mathcal{N}^{(\ell)} := \mathcal{A}^{(\ell)} - \mathcal{D}^{(\ell)}.
\end{equation*}
This leads to the Jacobi iteration
\begin{equation}
\label{eq:Jacobi}
x^{(\ell),n+1} = x^{(\ell),n} + \left( \mathcal{D}^{(\ell)} \right)^{-1} \odot \Big( b^{(\ell)} - \mathcal{A}^{(\ell)} x^{(\ell),n} \Big),
\end{equation}
where $\odot$ denotes the Hadamard product. Alternatively, we may use the damped form
\begin{equation}
\label{eq:Damped Jacobi}
x^{(\ell),n+1} = x^{(\ell),n} + \omega \left( \mathcal{D}^{(\ell)} \right)^{-1} \odot \Big( b^{(\ell)} - \mathcal{A}^{(\ell)} x^{(\ell),n} \Big), \quad 0 < \omega \leq 1.
\end{equation}
The inverse of $\mathcal{D}^{(\ell)}$ appearing in \eqref{eq:Jacobi} and \eqref{eq:Damped Jacobi} is interpreted elementwise; that is
\begin{equation*}
    \left( \mathcal{D}^{(\ell)} \right)^{-1}_{i_1,\ldots,i_d}
    := 1 / \mathcal{D}^{(\ell)}_{i_1,\ldots,i_d}.
\end{equation*}

In the context of low-rank methods, forming the elementwise inverse of the diagonal operator $\mathcal{D}^{(\ell)}$ is nontrivial, since this operation does not, in general, preserve low-rank tensor structure. However, several approximations can be employed. The simplest approach is to bound the spectrum of $\mathcal{D}^{(\ell)}$, yielding a Richardson-type iteration that ensures contractivity of the fixed-point iteration in \eqref{eq:Jacobi}, as in \cite{ballani2013projection}. While simple, this approach can sometimes lead to slow convergence. For the examples considered in this work, we find this simple Richardson approximation to be sufficient. When additional structure is available, more accurate approximations may be constructed. For example, if a constant-coefficient diffusion equation is discretized on a tensor-product grid, the resulting operator $\mathcal{A}^{(\ell)}$ admits a Kronecker-sum representation. In this case, the elementwise inverse is given by
\begin{equation*}
    \left( \mathcal{D}^{(\ell)} \right)^{-1} = \int_{0}^{\infty} e^{-t \mathcal{D}^{(\ell)}}\,dt = \int_{0}^{\infty} \bigotimes_{\mu=1}^{d} e^{-t D_{\mu}^{(\ell)}}\,dt,   
\end{equation*}
which can be approximated using a quadrature rule \cite{grasedyck2020parameter,grasedyck2023computing,grasedyck2025operator}. In the general case, one may instead use Newton--Schulz iteration, which defines this inverse as the root of the nonlinear equation
\begin{equation*}
    \mathcal{F}\left( \mathcal{X}_{i_1,\ldots,i_d}^{(\ell)} \right) = \mathcal{D}_{i_1,\ldots,i_d}^{(\ell)} - \frac{1}{\mathcal{X}_{i_1,\ldots,i_d}^{(\ell)}} = 0 \implies \mathcal{X}_{i_1,\ldots,i_d}^{(\ell),n+1} = \mathcal{X}_{i_1,\ldots,i_d}^{(\ell),n} \cdot \left( 2 - \mathcal{D}_{i_1,\ldots,i_d}^{(\ell)} \cdot \mathcal{X}_{i_1,\ldots,i_d}^{(\ell),n}  \right).
\end{equation*}
In the low-rank setting, these elementwise operations correspond to Hadamard products, which necessitate efficient algorithms to control rank growth \cite{kressner2014algorithm}. In recent work, Zheng et al. \cite{zheng2025HTACA} compute an elementwise inverse in the HT format using an adaptive cross-approximation technique with greedy sampling for a particular application involving spectral methods. A related strategy was employed in earlier work by Andreev and Tobler \cite{andreev2015multilevel} to approximate the action of an elementwise inverse on a residual within a smoothing procedure. Approaches based on cross approximation are particularly attractive for more general problems with nontrivial diagonals, as they require only pointwise access to tensor entries in order to construct the approximation.

At the coarsest level $\ell = 1$, the system may be solved using several possible approaches. A simple option is to expand the coarse-grid problem into its corresponding matrix form. Using the tensor product representation \eqref{eq:A as tensor product}, the system can be written as
\begin{equation*}
    \left( \sum_{k=1}^{r_A} \bigotimes_{\mu=1}^{d} A_{\mu,k}^{(1)} \right)
    \mathrm{vec}\!\left(x^{(1)}\right)
    =
    \mathrm{vec}\!\left(b^{(1)}\right),
\end{equation*}
where $\mathrm{vec}(\cdot)$ denotes the \emph{full} vectorization of the tensor arguments on the coarsest grid. This system may then be solved using either a direct or an iterative method, followed by reshaping and recompressing the solution into a low-rank format. While straightforward, this approach becomes impractical as the dimensionality or mode sizes increase. To preserve scalability and maintain the low-rank structure, we instead solve the coarse-grid system using the FGMRES method introduced in the previous section, without additional preconditioning. Since the coarsest problem is small by construction, this choice incurs only a modest computational cost while avoiding explicit formation of the full operator. We note that other low-rank linear solvers could instead be used, such as the AMEn method \cite{dolgov2014amen}, but we leave this investigation to future work.

A multigrid solver is defined through a recursive application of smoothing, restriction, coarse-grid correction, and prolongation, which characterize a \emph{cycle}. By varying the number and ordering of recursive coarse-grid corrections, one can control how computational effort is distributed across the hierarchy. A commonly used example is the V-cycle, which performs a single recursive coarse-grid correction per level and is typically the most efficient choice in terms of computational cost. The pseudo-code for a single V-cycle consists of the following steps:
\begin{enumerate}
    \item Pre-smoothing: apply $\nu_1$ iterations of the smoother, e.g., \eqref{eq:Damped Jacobi}.
    \item Compute the residual $r^{(\ell)} = b^{(\ell)} - \mathcal{A}^{(\ell)} x^{(\ell)}$.
    \item Restrict the residual: $r^{(\ell-1)} = \mathcal{R}^{(\ell)} r^{(\ell)}$.
    \item Coarse-grid correction: approximately solve $$ \mathcal{A}^{(\ell-1)} e^{(\ell-1)} = r^{(\ell-1)} $$ using a recursive multigrid cycle of the same type. At the coarsest level, the system is solved using the FGMRES method described in the previous section.
    \item Prolongate and correct: $$ x^{(\ell)} \leftarrow x^{(\ell)} + \mathcal{P}^{(\ell-1)} e^{(\ell-1)}. $$
    \item Post-smoothing: apply $\nu_2$ iterations of the smoother, e.g., \eqref{eq:Damped Jacobi}.
\end{enumerate}
More aggressive cycles can be obtained by increasing the number of recursive coarse-grid corrections. In particular, a W-cycle performs two recursive corrections per level, which improves robustness at the expense of additional computational work. An F-cycle provides an intermediate strategy, combining these approaches by performing an initial recursive correction followed by a V-cycle on the current level. In this work, we primarily consider preconditioners based on the V-cycle, while F- and W-cycles are used for comparison as standalone solvers. Further details on multigrid cycle implementations can be found in Briggs et al. \cite{briggs2000multigrid} and in Trottenberg et al. \cite{trottenberg2000multigrid}.

\section{Iterative Methods for Nonlinear Systems}
\label{sec:nonlinear solvers}

In this section we provide details regarding the iterative methods used to solve high-dimensional nonlinear systems of equations. We begin with a brief introduction of inexact Newton methods in the context of full-grid methods. Afterwards, we describe extensions of this approach which enable these methods to work in the low-rank setting. When describing these methods, we connect them to the linear solvers and preconditioning strategy introduced in the previous section.

\subsection{Inexact Newton Methods}
\label{subsec:inexact Newton methods}


Consider the following nonlinear system
\begin{equation}
    \label{eq:nonlinear tensor system}
    \mathcal{F}(u) = 0,
\end{equation}
where $u \in \mathbb{C}^{\mathbb{I}}$ and $\mathcal{F}$ denotes the nonlinear residual. Such systems arise from implicit discretizations of nonlinear PDEs and are often treated using Newton's method. In Newton's method, we start from an initial guess $u^{(0)}$ and generate a sequence of solutions to linearized systems
\begin{equation}
    \label{eq:linear tensor system for Newton}
    \mathcal{J}(u^{(k)}) \, \delta u^{(k)} = -\mathcal{F}(u^{(k)}), 
\end{equation}
which are then used to update the solution according to
\begin{equation}
    \label{eq:Newton correction equation}
    u^{(k+1)} = u^{(k)} + \delta u^{(k)}.
\end{equation}
Here, $\mathcal{J}(u^{(k)}) \in \mathbb{C}^{\mathbb{I} \times \mathbb{I}}$ denotes the Jacobian of $\mathcal{F}$ which is evaluated at $u^{(k)}$. In practice, the Newton iteration is terminated when the nonlinear residual satisfies either prescribed absolute or relative tolerances
\begin{equation}
    \label{eq:Newton stopping residual tolerances}
    \left\|\mathcal{F}\left( u^{(k)} \right)\right\|_F \leq \tau_{\text{abs}}, \quad \left\|\mathcal{F}\left( u^{(k)} \right)\right\|_F \leq \tau_{\text{rel}} \left\|\mathcal{F}\left( u^{(0)} \right)\right\|_F,
\end{equation}
where $\left\|\cdot \right\|_F$ is the Frobenius norm, or when the Newton updates become sufficiently small
\begin{equation}
    \label{eq:Newton stopping update tolerances}
    \left\| \delta u^{(k)} \right\|_F \leq \xi_{\text{abs}}, \quad \left\| \delta u^{(k)} \right\|_F \leq \xi_{\text{rel}} \left\| u^{(k)} \right\|_F.
\end{equation}

If the Jacobian is inexpensive to form, the linear systems \eqref{eq:linear tensor system for Newton} can be treated directly. However, in high-dimensional problems, forming the Jacobian may be
computationally prohibitive or infeasible in problems with coupled nonlinearities. In such situations, one can instead use Jacobian-free Newton--Krylov (JFNK) methods \cite{knoll2004jfnk}, which approximately solve the linear systems \eqref{eq:linear tensor system for Newton} in each Newton iteration using a Krylov method such as FGMRES \cite{saad1993flexible}. JFNK methods avoid forming the Jacobian explicitly and instead require only Jacobian-vector products of the form $\mathcal{J}(u^{(k)})v$. These products may be supplied analytically when the action of the Jacobian is available, or approximated through a finite-difference perturbation of the nonlinear residual \eqref{eq:nonlinear tensor system}. In the latter case, one uses
\begin{equation}
    \label{eq:Jacobian-vector product approximation}
    \mathcal{J}(u)v 
    \approx 
    \frac{\mathcal{F}(u + \varepsilon v) - \mathcal{F}(u)}{\varepsilon},
\end{equation}
where $\varepsilon$ defines a small perturbation of $u$ that must be carefully chosen to balance approximation errors and avoid underflow issues. A comprehensive overview of choices for this parameter is given in
\cite{knoll2004jfnk}. One common approach selects $\varepsilon$ dynamically based on the sizes of $u$ and $v$:
\begin{equation*}
    \varepsilon = \frac{ \sqrt{(1 + \left\| u \right\|_F)\epsilon_{\text{mach}}} }{ \left\| v \right\|_F }.
\end{equation*}
Here, $\epsilon_{\text{mach}}$ is the machine precision.

The inner linear solver in this work is the FGMRES method described in Section~\ref{subsec:fgmres} with the GMG method described in Section~\ref{subsec:multigrid_ht} as a preconditioner. As the Jacobian depends on the current Newton state, the user may wish to update the preconditioner at the beginning of each Newton iteration. In order to avoid over-solving the inner linear systems \eqref{eq:linear tensor system for Newton} during the early Newton iterations, we adapt the relative tolerance according to
\begin{equation}
    \label{eq:Inner Newton forced relative tolerance}
    \left\| \mathcal{J}(u^{(k)}) \, \delta u^{(k)} + \mathcal{F}(u^{(k)}) \right\|_F \leq \gamma_{k} \left\| \mathcal{F}(u^{(k)}) \right\|_F .
\end{equation}
The parameter $\gamma_k \in (0,1)$ is a forcing term that controls the accuracy of the linear solve. In this work, we employ the heuristic
\begin{equation}
    \label{eq:forcing term}
    \gamma_k =
    \min\left(\gamma_{\max},~
    \max\left(\gamma_{\min},~
    \frac{1}{2}\sqrt{\frac{\left\| \mathcal{F}(u^{(k)}) \right\|_F}{\left\| \mathcal{F}(u^{(0)}) \right\|_F}}
    \right)\right),
\end{equation}
where $\gamma_{\min}$ and $\gamma_{\max}$ are prescribed bounds on the forcing term. We note that other choices can be used for inexact Newton methods \cite{eisenstat1996choosing}.


The solution of the linear systems \eqref{eq:linear tensor system for Newton}, and the corresponding update \eqref{eq:Newton correction equation}, may fail to produce a sufficient decrease in the nonlinear residual. This can occur when the initial guess is far from the solution or when the linear systems are solved inexactly during the early Newton iterations. To improve robustness, we employ a line-search strategy in which the update \eqref{eq:Newton correction equation} is modified according to
\begin{equation}
    \label{eq:line search for Newton}
    u^{(k+1)} = u^{(k)} + \alpha_{k}\delta u^{(k)},
\end{equation}
where the step length $\alpha_{k} \in (0,1]$ is chosen to ensure sufficient decrease in the nonlinear residual. Specifically, $\alpha_{k}$ is determined via a backtracking procedure that enforces an Armijo-type condition
\begin{equation}
    \label{eq:Armijo backtracking}
    \left\| \mathcal{F}(u^{(k)} + \alpha_{k} \delta u^{(k)}) \right\|_F
    \leq
    (1 - c\alpha_{k}) \left\| \mathcal{F}(u^{(k)}) \right\|_F,
\end{equation}
for some constant $c \in (0,1)$. A simple strategy for updating $\alpha_{k}$ in the backtracking step is to initially define $\alpha_{k}=1$ and recursively halve its value until the condition \eqref{eq:Armijo backtracking} is satisfied. The Armijo condition is satisfied for a sufficiently small $\alpha_{k}$ under the standard assumption that the computed Newton correction is a descent direction for the nonlinear residual. In our implementation, we impose the lower bound $\alpha_{k} \geq 10^{-4}$, which guarantees that the backtracking procedure terminates after finitely many reductions.

\subsection{Inexact Newton Methods for Tensor Decompositions}
\label{subsec:inexact Newton tensor methods}


As discussed in Section~\ref{subsec:fgmres tensor}, we seek to represent the solution and all intermediate quantities in a compressed tensor format. The only requirement on the representation is that it supports basic operations such as addition, scalar multiplication, truncation, inner products, and norm evaluation. Within this framework, the inexact Newton method introduces several additional considerations. In particular, the evaluation of the nonlinear residual $\mathcal{F}(u)$ and the Jacobian-vector products required by the linear solver must be carried out in the presence of truncation. As a result, these operations are no longer exact, and the computed quantities are subject to perturbations introduced by low-rank approximation. 

One such modification arises in the Jacobian-vector product approximation \eqref{eq:Jacobian-vector product approximation}, which can be approximated as
\begin{equation*}
    \mathcal{J}(u)v \approx  \mathfrak{T}\left( \frac{\mathcal{F}(u + \varepsilon v) - \mathcal{F}(u) }{\varepsilon}\right).
\end{equation*}
That is, the residuals are evaluated at the untruncated states, and truncation is applied only to their difference. This approach avoids the possibility of prematurely truncating the perturbation, which might lead to ill-conditioning and degrade the approximation of the Jacobian-vector product. These perturbations carry over to the inner linear solver, which, as described in Section~\ref{subsec:fgmres tensor}, already departs from the standard linear theory due to the use of truncation. As a result, the linear systems arising in each Newton step are solved only approximately.

An important consequence is the interaction between truncation and the inexact Newton forcing term \eqref{eq:forcing term}. The forcing term $\gamma_k$ controls the accuracy of the linear solve via \eqref{eq:Inner Newton forced relative tolerance}, while truncation imposes a lower bound on the attainable accuracy. Solving the linear systems to a tolerance significantly smaller than the truncation error does not improve the quality of the Newton update. Therefore, the forcing term and truncation tolerances must be chosen in a consistent manner. Moreover, truncation implies that the nonlinear residual is evaluated only approximately, and thus the line-search procedure may be influenced by truncation errors, particularly when the residual norm is small. Similar to the approximation of Jacobian-vector products, we modify the backtracking procedure \eqref{eq:line search for Newton} to incorporate truncation:
\begin{equation*}
    u^{(k+1)} = \mathfrak{T} \Big( u^{(k)} + \alpha_{k}\delta u^{(k)} \Big).
\end{equation*}

Lastly, users should select the truncation tolerances so that they are compatible with the target accuracy of the nonlinear solve. An additional practical consideration is the numerical stability of tensor operations involving quantities of small magnitude. In particular, when the residual $\mathcal{F}(u^{(k)})$ becomes small, the accuracy of norm and inner product evaluations may degrade if the tensor representation is not suitably conditioned \cite{kressner2014algorithm}. To mitigate this issue, it is advantageous to apply orthogonalization procedures to the tensor factors \emph{prior} to performing such operations. This improves the numerical stability of the representation and ensures that quantities such as residual norms and inner products are computed with sufficient accuracy. In practice, this step is especially important in the later stages of the Newton iteration, where small inaccuracies can adversely affect both the stopping criteria and the behavior of the line search.

\section{Numerical Examples}
\label{sec:numerics}

In this section, we apply the proposed iterative methods to several well-known test problems from the literature to study their convergence properties. The examples are arranged in increasing order of complexity and aim to benchmark the capabilities of the proposed algorithms. All methods were implemented in MATLAB using the HTucker package \cite{kressner2014algorithm}, which provides the underlying data structures for the HTD, as well as methods for tensor arithmetic, orthogonalization, truncation, and contractions. Unless stated otherwise, we use $\epsilon_{\mathrm{abs}} = \epsilon_{\mathrm{rel}} = 10^{-4}$ for the absolute and relative truncation tolerances in all examples. Recall that these tolerances are used during truncation to select the hierarchical ranks according to \eqref{eq:abs trunc tolerance} and \eqref{eq:rel trunc tolerance}, respectively. In addition, all experiments set $\eta_{\mathrm{rel}} = 0.1$ in the GMRES/FGMRES stopping criterion \eqref{eq:FGMRES inner acceptance criterion}. Additionally, the adaptive-rank GMG components use the iterative method from Section \ref{subsec:fgmres tensor} without preconditioning on the coarsest level, and the number of levels is chosen so that the coarsest grid has mode size $9$.

\subsection{Poisson's Equation}
\label{subsec:poisson}

As a first example, we consider the $d$-dimensional Poisson problem on the unit cube $\Omega_{\mathbf{x}} = [0,1]^d$ with homogeneous Dirichlet boundary conditions:
\begin{equation}
\label{eq:Poisson problem}
\left\{
\begin{aligned}
    -\Delta_{\mathbf{x}} u &= f(\mathbf{x}), &&\mathbf{x} \in \Omega_{\mathbf{x}}, \\
    u(\mathbf{x}) &= 0, &&\mathbf{x} \in \partial\Omega_{\mathbf{x}}.
\end{aligned}
\right.
\end{equation}
Poisson's equation appears in a variety of applications and is frequently used to study the convergence properties of iterative methods.

To assess accuracy and convergence, we employ a manufactured solution consisting of a smooth mixture of Laplacian eigenmodes:
\begin{equation*}
    u(\mathbf{x})
    = \sum_{m=1}^{M}c_{m}\prod_{\mu=1}^d \sin(k_{m,\mu}\pi x_\mu),
\end{equation*}
where
\begin{equation*}
    k_{m,\mu} \sim \mathcal{U}(1,k_{\text{max}}), \quad c_{m} \sim \mathcal{N}(0,1),
\end{equation*}
with $M = 5$ and $k_{\text{max}} = 4$. The corresponding right-hand side in \eqref{eq:Poisson problem} is
\begin{equation*}
    f(\mathbf{x})
    = \sum_{m=1}^{M}c_{m}\lambda_{m}\prod_{\mu=1}^d \sin(k_{m,\mu}\pi x_\mu), \qquad \lambda_{m} := \pi^{2}\sum_{\mu=1}^{d}k_{m,\mu}^{2}.
\end{equation*}
This construction satisfies the boundary conditions and avoids alignment with a single coordinate direction, preventing trivial convergence of iterative methods.

The domain is discretized using a uniform, vertex-centered grid with $N = 2^{k} + 1$ points in each dimension which provides a proper nesting of the grids across multigrid levels. The Laplacian is approximated using second-order central finite differences
\begin{equation*}
    \Delta_{\mathbf{x}}u \approx \sum_{\mu = 1}^{d} D_{x_{\mu}x_{\mu}}u,
\end{equation*}
which are applied in a matrix-free manner along the corresponding modes of the tensor.

\begin{figure}[!htbp]
\centering

\begin{subfigure}{0.325\linewidth}
\centering
\includegraphics[width=\linewidth]{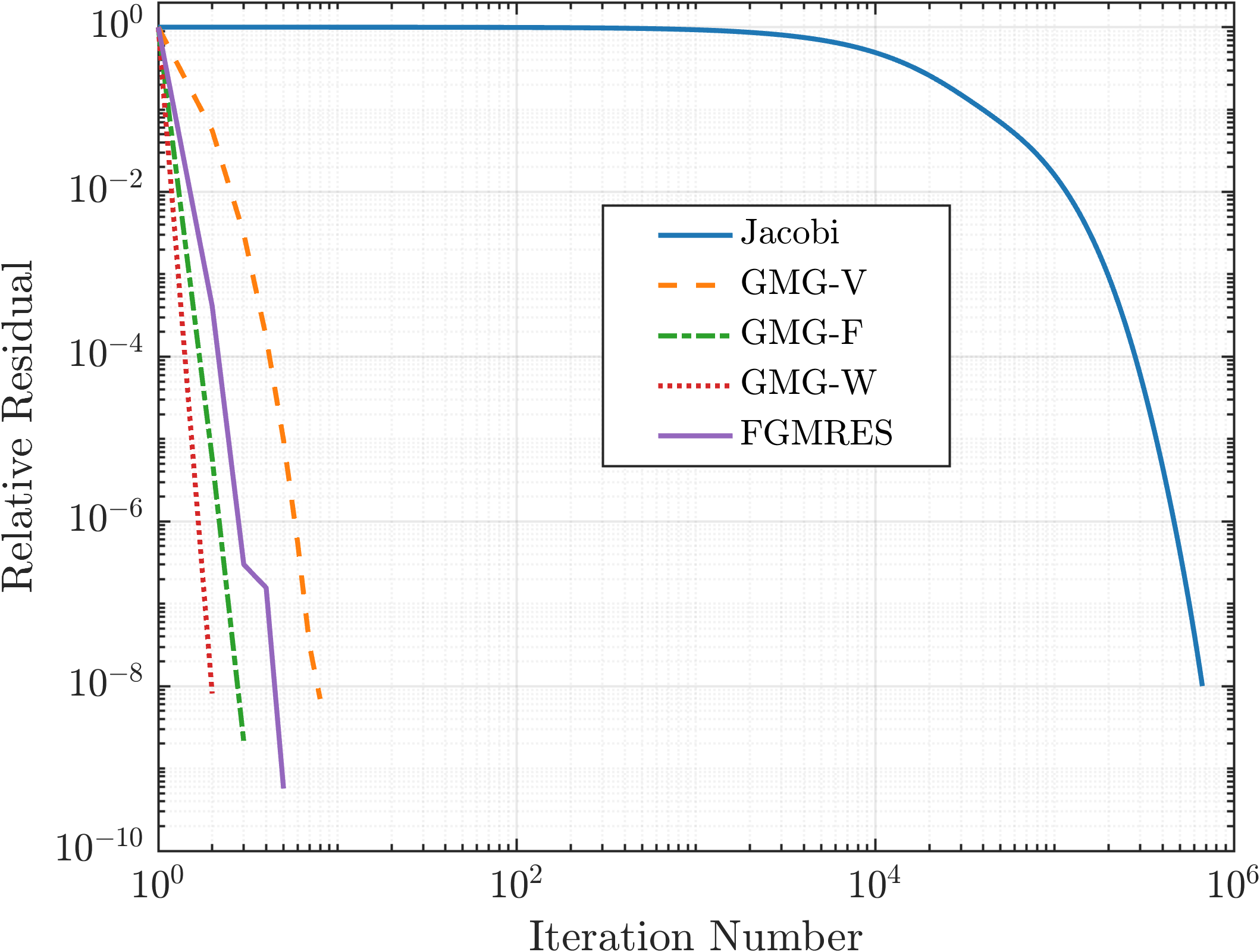}
\caption{Solver comparison}
\label{fig:poisson solver comparison}
\end{subfigure}
\begin{subfigure}{0.32\linewidth}
\centering
\includegraphics[width=\linewidth]{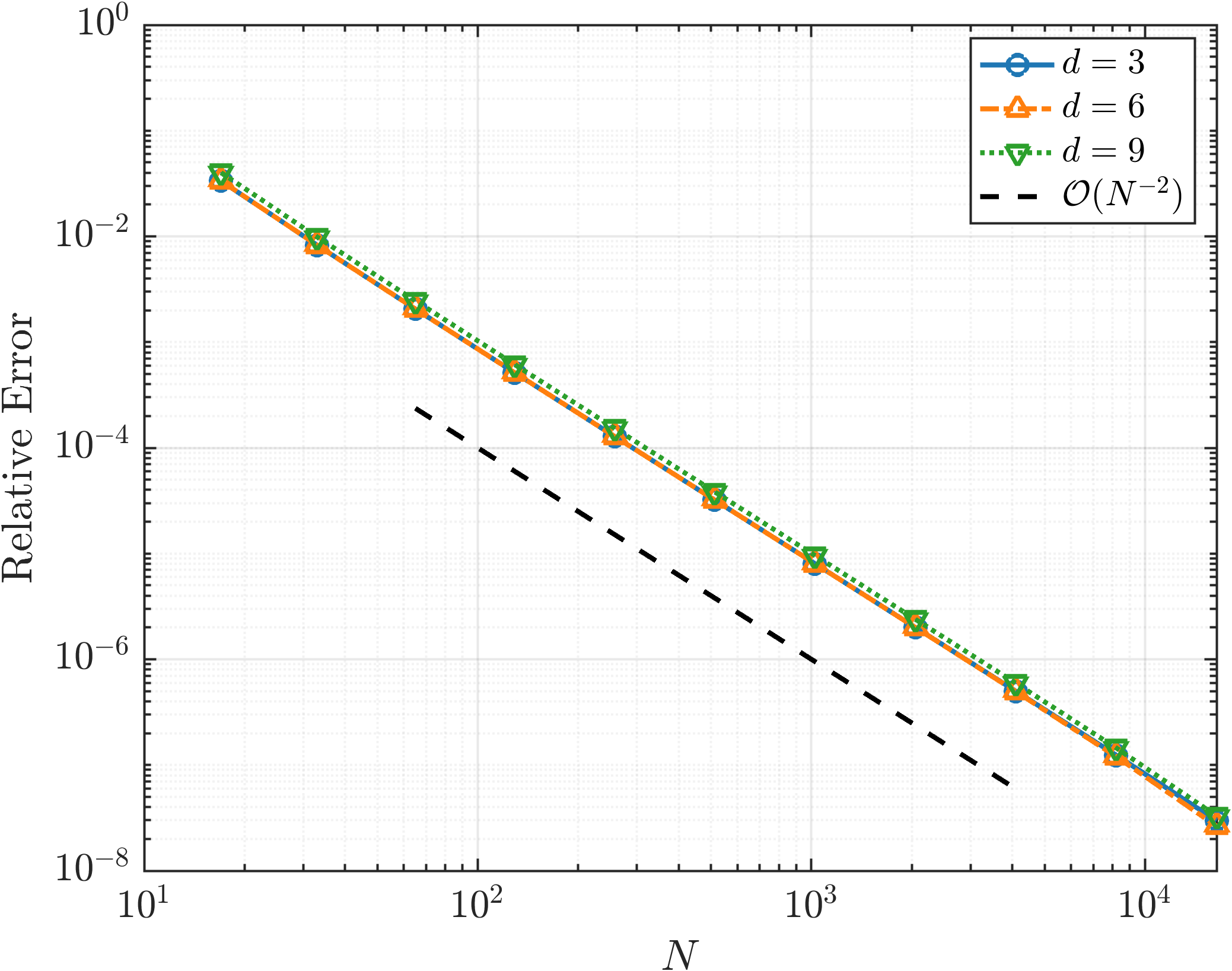}
\caption{GMG-V mesh refinement}
\label{fig:poisson mesh refinement}
\end{subfigure}
\begin{subfigure}{0.32\linewidth}
\centering
\includegraphics[width=\linewidth]{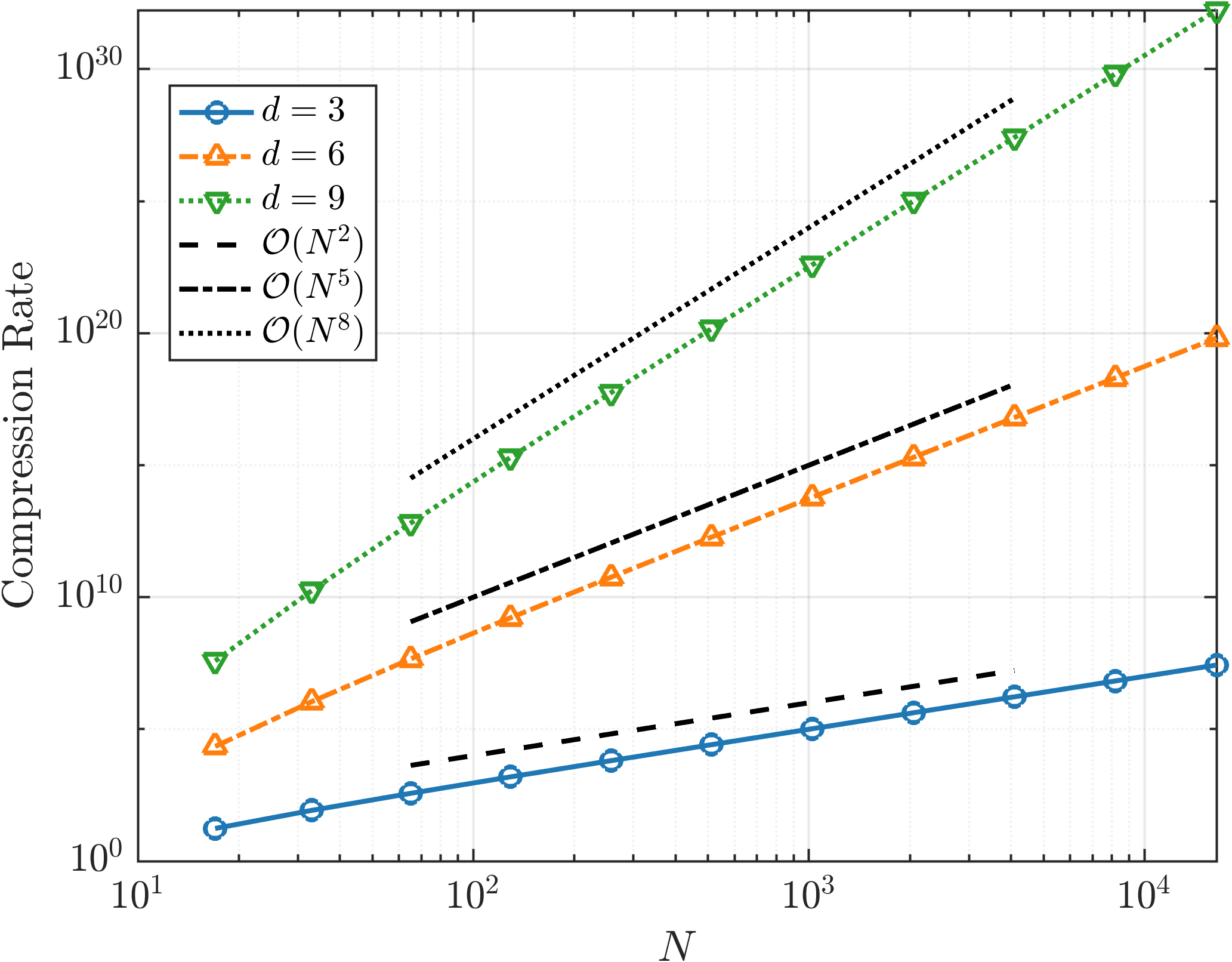}
\caption{GMG-V compression rate}
\label{fig:poisson compression}
\end{subfigure}

\begin{subfigure}{0.325\linewidth}
\centering
\includegraphics[width=\linewidth]{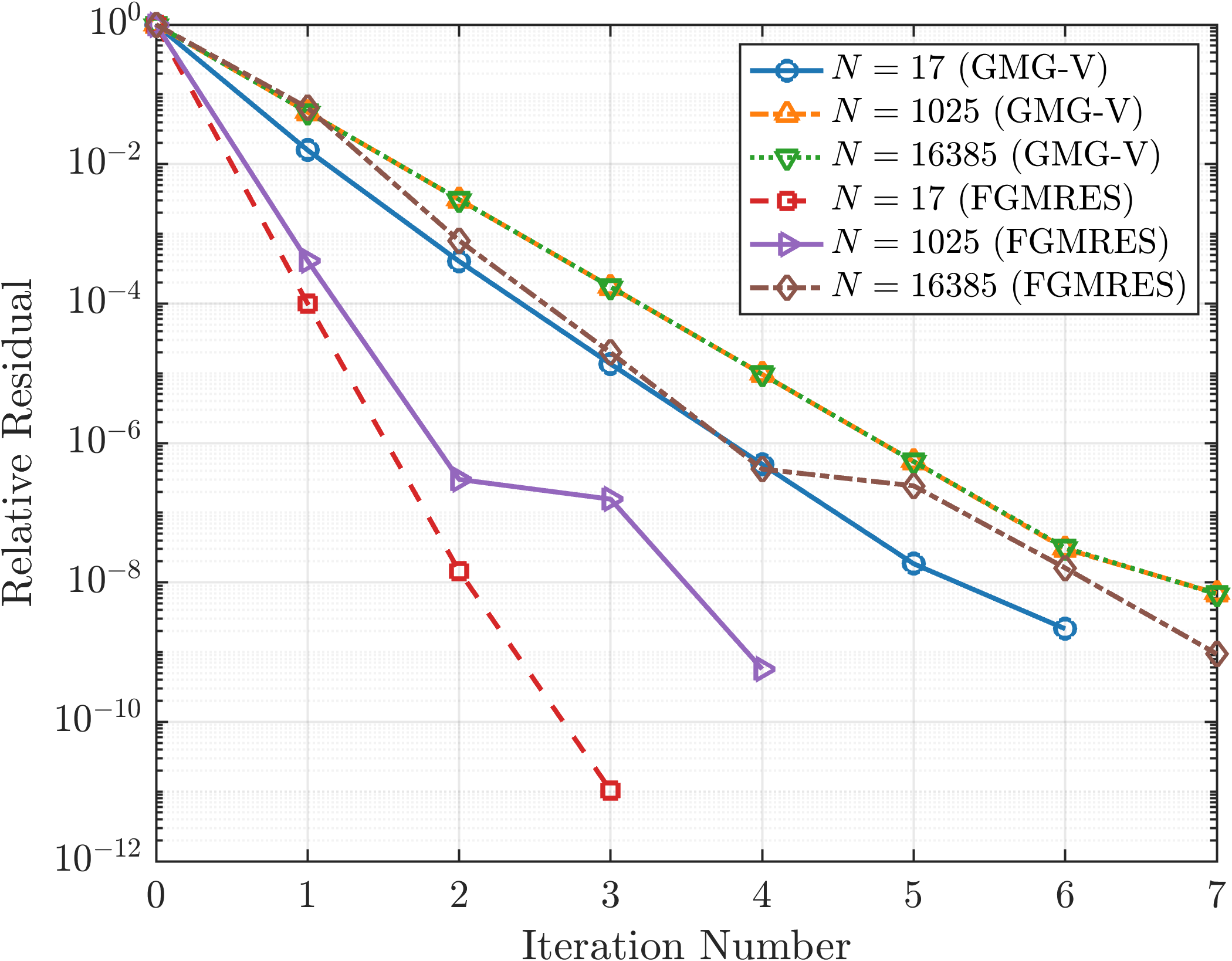}
\caption{$N$ dependence ($d=3$)}
\label{fig:poisson mesh independence}
\end{subfigure}
\begin{subfigure}{0.325\linewidth}
\centering
\includegraphics[width=\linewidth]{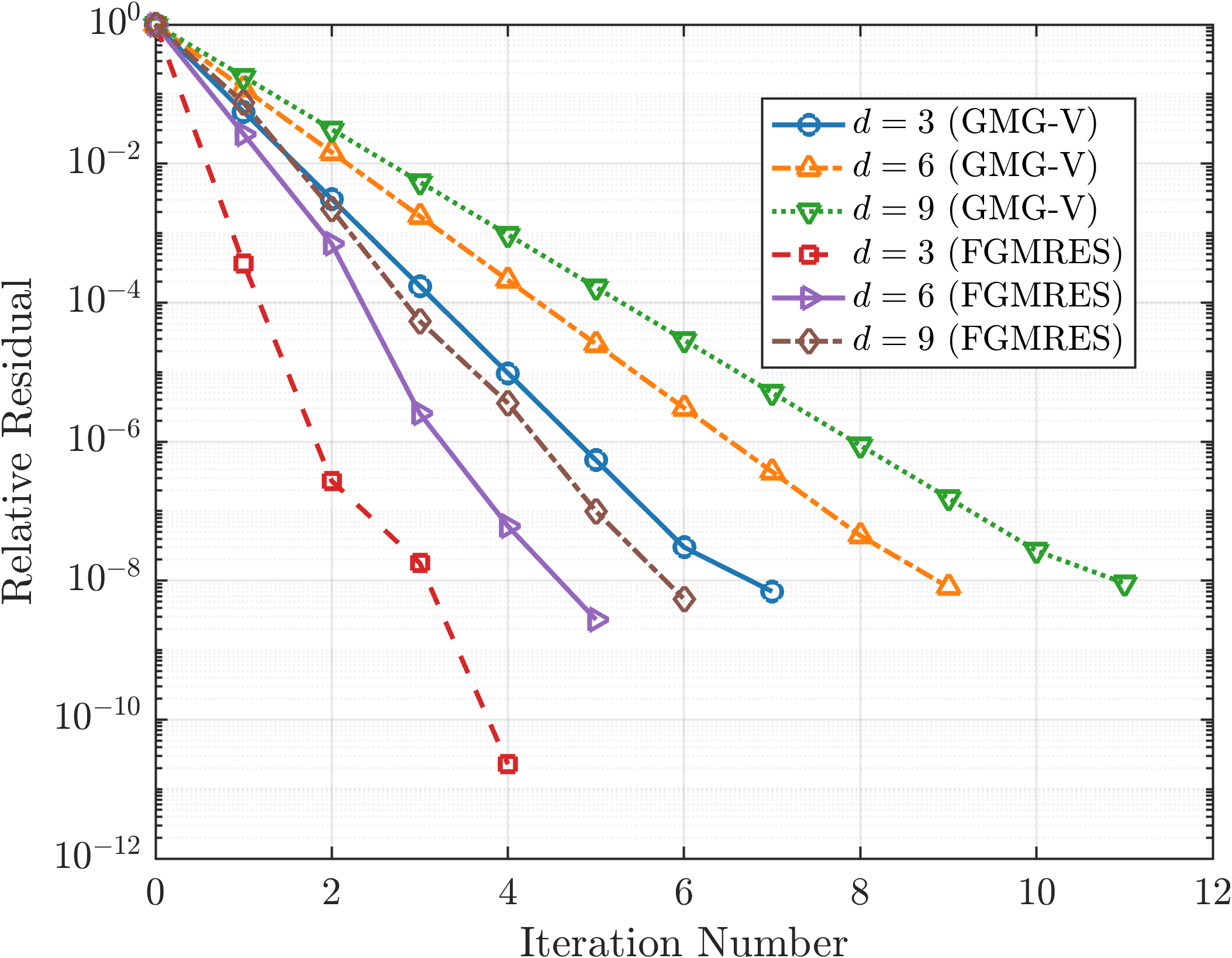}
\caption{$d$ dependence ($N=513$)}
\label{fig:poisson dim independence}
\end{subfigure}
\begin{subfigure}{0.3175\linewidth}
\centering
\includegraphics[width=\linewidth]{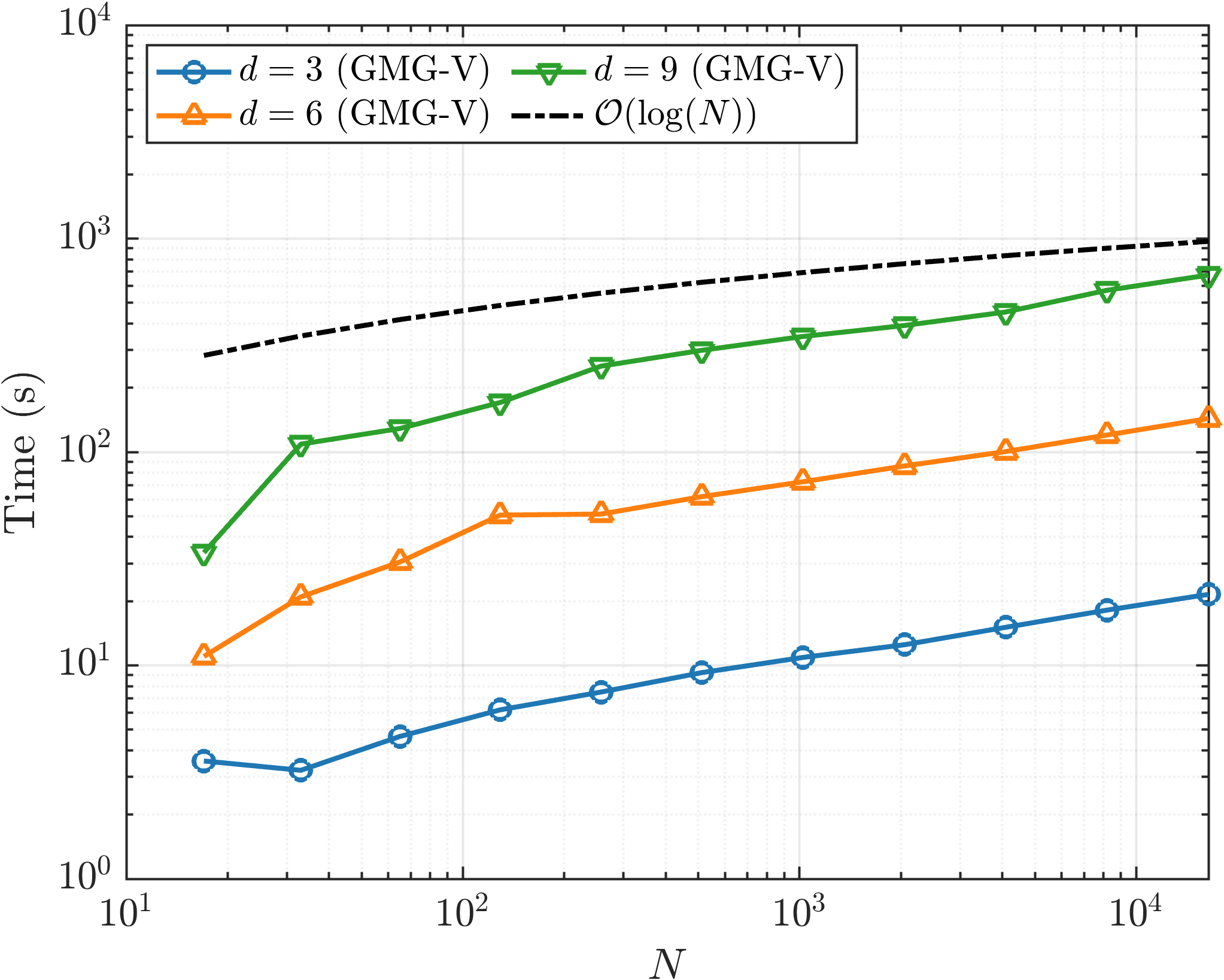}
\caption{GMG-V complexity}
\label{fig:poisson solver complexity}
\end{subfigure}

\caption{(Poisson problem)
Comparison of adaptive-rank solvers and scalability properties.
(a) Relative residual versus iteration for Jacobi, GMG (V/F/W cycles), and GMG-V preconditioned FGMRES on a fixed mesh ($N=1025$, $d=3$), demonstrating the rapid convergence of multilevel methods.
(b) Mesh refinement study for GMG-V methods confirming second-order accuracy of the finite-difference discretization. Results are shown only for the GMG-V, as other methods produced identical results.
(c) Compression rate versus $N$ for several values of $d$ for the GMG-V method, illustrating the reduction in memory achieved by the adaptive-rank representation. We only show results for the GMG-V method, as other methods produced identical results.
(d) Relative residual histories for GMG-V and GMG-V preconditioned FGMRES using a fixed dimension $d = 3$ with different per-dimension mode sizes $N$. The iteration count remains essentially independent of the per-dimension mode size.
(e) Relative residual histories for GMG-V and GMG-V preconditioned FGMRES using a fixed per-dimension mode size $N = 513$ with different dimensions $d$. The iteration count shows a weak dependence on the dimension.
(f) The solver time for the GMG-V method as a function of the per-dimension mode size $N$ for different dimensions $d$. Over the tested range, the measured growth is well described by a reference curve proportional to $\mathcal{O}(\log (N))$, which is more favorable than the conservative $\mathcal{O}(N\log (N))$ estimate suggested by the mesh hierarchy alone. This behavior is likely influenced by the adaptive-rank representation and by implementation-dependent effects, with the cost of traversing the multilevel hierarchy dominating the observed growth.
}
\label{fig:poisson all results}
\end{figure}

Figure \ref{fig:poisson solver comparison} compares the performance of several adaptive-rank linear solvers for the Poisson problem using a fixed dimension $d =3$ and per-dimension mode size $N = 1025$. After each inter-grid transfer, we apply 10 passes of the undamped Jacobi method as a smoother. On each level $\ell$, we approximate the inverse of the diagonal operator by a scalar multiple of the identity:
\begin{equation*}
    \left(\mathcal{D}^{(\ell)}\right)^{-1} \approx \left[ 
     2\sum_{\mu=1}^d \frac{1}{\left(h_{\mu}^{(\ell)}\right)^2} \right]^{-1} \mathcal{I}.
\end{equation*} 
We compare the relative residual history obtained with Jacobi, GMG using V-, F-, and W-cycles, as well as an FGMRES solver preconditioned with the GMG-V method. Each method converges the solution until the relative residual error is below $10^{-8}$. Among these methods, the proposed GMG solvers exhibit rapid convergence, which suggests that classical GMG behavior holds in the adaptive-rank setting as well. In particular, we find that 1 W-cycle is sufficient to converge the solution, while the F-cycle and V-cycle use 2 and 7 iterations, respectively. Note that these iteration counts ignore the iterations associated with the solve on the coarsest grid, as these problems are quite small and rapidly converge within the same tolerances used for the outer solver.

The results of a standard mesh refinement study are provided in Figure \ref{fig:poisson mesh refinement}. In this experiment the per-dimension mode size is set according to $N = 2^{k} + 1$, and we set $k = 4,5, \cdots,14$. This creates a collection of problems with mode sizes ranging from $N = 17$ to $N = 16,385$. The study was also conducted with several different dimensions, namely $d = 3, 6, 9$ and confirms the second-order accuracy of the finite-difference approximations. For brevity, only results for the GMG solver using a V-cycle are shown. Similar results are obtained with the F- and W-cycle, so we omit them.

The compression behavior of the adaptive-rank representation using the GMG method with a V-cycle is shown in Figure \ref{fig:poisson compression}, which demonstrates the benefit of the adaptive-rank discretization. This plot shows the compression rate, provided in Definition \ref{def:compression_rate}, as a function of the per-dimension mode size $N$, under several different choices for the dimension $d$. The plot also includes several reference lines corresponding to the asymptotic compression rate predicted by the HT format discussed in Remark \ref{rem:asymptotic_compression}. Similar results are observed using the F- and W-cycle as well as preconditioned FGMRES method, so we omit them.

GMG methods often exhibit nearly mesh-independent convergence properties, meaning that the total number of cycles used depends weakly on the mesh resolution $N$. In addition, previous work by Ballani and Grasedyck \cite{ballani2013projection} suggests that the number of cycles also depends \emph{weakly} on the dimension $d$ of the problem. We present the results of a similar study in this work for the GMG-V solver as well as the FGMRES solver which uses the GMG-V method as a preconditioner in Figures \ref{fig:poisson mesh independence} and \ref{fig:poisson dim independence}. Figure \ref{fig:poisson mesh independence} presents the relative residual data collected in a mesh refinement study using a  fixed dimension $d = 3$. The results suggest that the number of iterations is independent of the per-dimension mode size, even for highly resolved meshes. Figure \ref{fig:poisson dim independence} presents the relative residual data obtained by fixing the per-dimension mode size $N = 513$ and varying the number of dimensions. We find that the total number of iterations used depends quite weakly on the dimensionality $d$. These results suggest that multilevel adaptive-rank methods can be effective at treating problems with high-resolution and dimensionality requirements.

Figure~\ref{fig:poisson solver complexity} presents the measured computational complexity of the GMG-V solver. For each dimension $d$, the total time to solution is plotted as a function of the per-dimension mode size $N$, along with a reference line corresponding to $\mathcal{O}(\log (N))$. The results indicate that the pure GMG solver exhibits approximately logarithmic complexity with respect to $N$ and a weak dependence on the dimensionality $d$. Together with the mesh- and dimension-independence results in Figures~\ref{fig:poisson mesh independence} and~\ref{fig:poisson dim independence}, these results indicate that the GMG-V solver remains effective for large per-dimension mode sizes. Based on the mesh hierarchy alone, a conservative estimate would suggest a cost of order $\mathcal{O}(N\log (N))$ with respect to the per-dimension mode size. The more favorable empirical behavior observed here is likely influenced by the adaptive-rank representation and by implementation-dependent effects, including MATLAB overheads and the relatively low cost of one-dimensional operations used in the test. We note that similar behavior has been observed in other high-dimensional low-rank solvers (see, e.g., \cite{sands2025transport,zheng2025HTACA}). In this regime, the cost associated with traversing the multilevel hierarchy appears to dominate the measured growth, leading to an approximately logarithmic dependence on $N$.

\subsection{Dougherty--Fokker--Planck Equation}
\label{subsec:dougherty}

The Dougherty--Fokker--Planck (DFP) equation \cite{dougherty1964model} is a nonlinear convection-diffusion equation which can be used to model collisional dynamics of particles found in both neutral gases and plasmas. We consider a spatially homogeneous (0D–3V) problem on the truncated velocity domain $\Omega_{\mathbf{v}} = [-6,6]^{3}$ given by
\begin{equation}
\label{eq:DFP system}
\left\{
\begin{aligned}
    &\partial_{t}f = \nabla_{\mathbf{v}} \cdot \Big( T \nabla_{\mathbf{v}} f + (\mathbf{v} - \mathbf{u})f \Big),  &&(\mathbf{v},t) \in \Omega_{\mathbf{v}} \times (0,T_{f}], \\
    &f(\mathbf{v},0) = f_{0}(\mathbf{v}),  &&\mathbf{v} \in \Omega_{\mathbf{v}},
\end{aligned}
\right.
\end{equation}
subject to periodic boundary conditions. The temperature $T > 0$ and mean velocity $\mathbf{u}$ are defined through the moments
\begin{equation*}
    n = \int_{\Omega_{\mathbf{v}}} f \,d\mathbf{v}, \quad n\mathbf{u} = \int_{\Omega_{\mathbf{v}}} \mathbf{v}f \,d\mathbf{v}, \quad E = \frac{1}{2} \int_{\Omega_{\mathbf{v}}} |\mathbf{v}|^{2}f \,d\mathbf{v},
\end{equation*}
corresponding to particle number density, momentum density, and energy density. The energy can be written alternatively as
\begin{equation*}
    E = \frac{1}{2} \left( n |\mathbf{u}|^{2} + 3nT\right),
\end{equation*}
where the temperature $T$ is defined as the central moment
\begin{equation*}
    T = \frac{1}{3n} \int_{\Omega_{\mathbf{v}}} |\mathbf{v} - \mathbf{u}|^{2}f \,d\mathbf{v}.
\end{equation*}
Under periodic boundary conditions, the collision operator preserves these moments, implying $\partial_t n = \partial_t (n\mathbf{u}) = \partial_t E = 0$, and hence $\partial_t T = 0$. Therefore, in the 0D–3V setting, the problem reduces to a linear convection--diffusion equation with coefficients determined by the initial data.

We consider the bimodal Maxwellian initial condition
\begin{equation}
    \label{eq:bimaxwellian initial data}
    f_{0}(\mathbf{v}) = \frac{n_{1}}{\left( 2\pi T_{1}\right)^{3/2}} \exp \left( -\frac{|\mathbf{v} - \mathbf{u}_{1}|^{2}}{2T_{1}} \right) + \frac{n_{2}}{\left( 2\pi T_{2}\right)^{3/2}} \exp \left( -\frac{|\mathbf{v} - \mathbf{u}_{2}|^{2}}{2T_{2}} \right),
\end{equation}
where
\begin{equation}
    \label{eq:bimaxwellian params}
    n_{1}=\frac{1}{2}, \quad
    n_{2}=\frac{1}{2}, \quad
    \mathbf{u}_1=(2,2,0)^{T}, \quad
    \mathbf{u}_2=(-2,-2,0)^{T}, \quad
    T_{1}=2, \quad
    T_{2}=2.
\end{equation}
This initial condition has hierarchical rank 2 and relaxes toward the equilibrium Maxwellian
\begin{equation}
    \label{eq:equilibrium distribution}
    f_{\infty}(\mathbf v) = \frac{ n }{(2\pi T)^{3/2}} \exp \left(-\frac{|\mathbf{v} - \mathbf{u}|^2}{2T}\right),
\end{equation}
which has rank 1. Using the configuration given in equation \eqref{eq:bimaxwellian params}, the parameters of the equilibrium distribution are determined to be $n = 1$, $\mathbf{u} = \mathbf{0}$, and $T = 14/3.$

For temporal discretization, we apply backward Euler to equation \eqref{eq:DFP system}, which yields
\begin{equation*}
    f^{n+1} - \Delta t\nabla_{\mathbf v} \cdot 
    \Big(T \nabla_{\mathbf{v}} f^{n+1} + (\mathbf{v} - \mathbf{u}) f^{n+1}\Big)
    = f^{n}.
\end{equation*}
Again, $T$ and $\mathbf{u}$ are computed analytically from the initial condition. The implicit treatment of both the convection and diffusion components is motivated by the following observation. For explicit discretizations on a uniform velocity grid with mesh spacing $h$, the convection term yields the usual CFL stability restriction of the form
\begin{equation*}
    \Delta t_{\text{exp}} \lesssim \frac{h}{a}, 
\end{equation*}
where $a = \max_{\mu} |v_{\mu} - u_{\mu}|$ denotes the maximum characteristic speed along dimension $v_{\mu}$. In contrast, the diffusion term yields the more restrictive parabolic stability condition of the form
\begin{equation*}
    \Delta t_{\text{exp}} \lesssim \frac{h^{2}}{T},
\end{equation*}
which becomes prohibitively restrictive as the mesh is refined.

The discretization of the collision operator uses periodic mode-wise finite difference approximations with second-order accuracy. In order to treat the convection terms, we define the convective flux $F_{\mu}(f) = (v_{\mu}-u_{\mu})f$ and define $\alpha_{\mu} = | v_{\mu}-u_{\mu} |$. Applying a Lax-Friedrichs flux splitting gives
\begin{equation*}
    F_{\mu}^{\pm}(f) = \frac{1}{2}\Big( F_{\mu}(f) \pm \alpha_{\mu} f \Big), \quad \mu = 1,2,3,
\end{equation*}
where $F_{\mu}^{+} \geq 0$ and $F_{\mu}^{-} \leq 0$. Then, we approximate the first-order derivatives as
\begin{equation*}
    \partial_{v_{\mu}}F_{\mu}(f)
    \approx 
    D_{v_{\mu}}^{-} \big(F_{\mu}^{+}\big) + D_{v_{\mu}}^{+}\big(F_{\mu}^{-}\big), \quad \mu = 1,2,3,
\end{equation*}
where $D_{v_{\mu}}^{-}$ and $D_{v_{\mu}}^{+}$ represent second-order backward and forward difference approximations to $\partial_{v_{\mu}}$, respectively. The fully discrete DFP collision operator is
\begin{equation*}
    \mathcal{L} f
    =
    - T \sum_{\mu=1}^{3} D_{v_{\mu}v_{\mu}} f
    -
    \sum_{\mu=1}^{3}
    \left(
    D_{v_{\mu}}^{-} \big( F_{\mu}^{+}(f) \big)
    +
    D_{v_{\mu}}^{+} \big( F_{\mu}^{-}(f) \big)
    \right),
\end{equation*}
and the corresponding linear system to be solved at each time step is given by
\begin{equation}
    \left( \mathcal{I} + \Delta t \mathcal{L} \right) f^{n+1} = f^{n}.
    \label{eq:DFP_system_final}
\end{equation}

The linear system \eqref{eq:DFP_system_final} is solved using the proposed adaptive-rank FGMRES approach with GMG-V preconditioning. The preconditioner consists of 2 iterations of a V-cycle and is smoothed using 10 iterations of a damped Jacobi method, with $\omega = 0.7$. On level $\ell$, we approximate the elementwise inverse of the diagonal operator by a scalar multiple of the identity using upper bounds for the convective diagonal contributions. This gives
\begin{equation*}
    \left(\mathcal{D}^{(\ell)}\right)^{-1} \approx \left[ 1+ \Delta t
     \left(2T\sum_{\mu=1}^d \frac{1}{\left(h_{\mu}^{(\ell)}\right)^2} + \frac{3}{2}\sum_{\mu=1}^d \frac{\max_{\nu\in v_{\mu}^{(\ell)}}|\nu - u_{\mu}|}{h_{\mu}^{(\ell)}} \right) \right]^{-1} \mathcal{I}.
\end{equation*}
The methods converge the solution at each time step until the relative residual error is below $10^{-4}$.

\begin{figure}[!htbp]
    \centering
    \begin{subfigure}{0.32\linewidth}
        \centering
        \includegraphics[width=\linewidth]{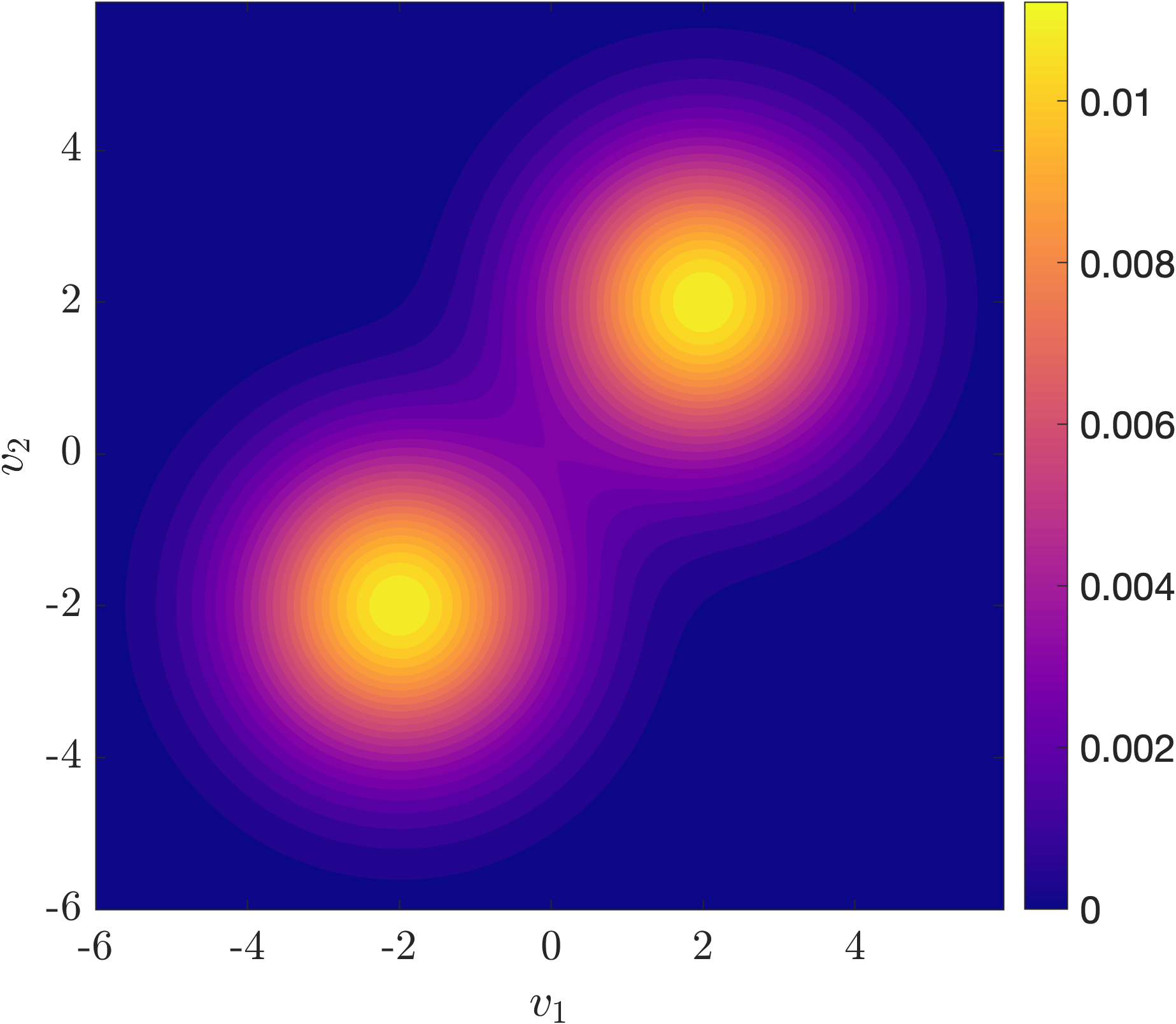}
        \caption{$t = 0$}
        \label{fig:dfp t = 0}
    \end{subfigure}
    \begin{subfigure}{0.32\linewidth}
        \centering
        \includegraphics[width=\linewidth]{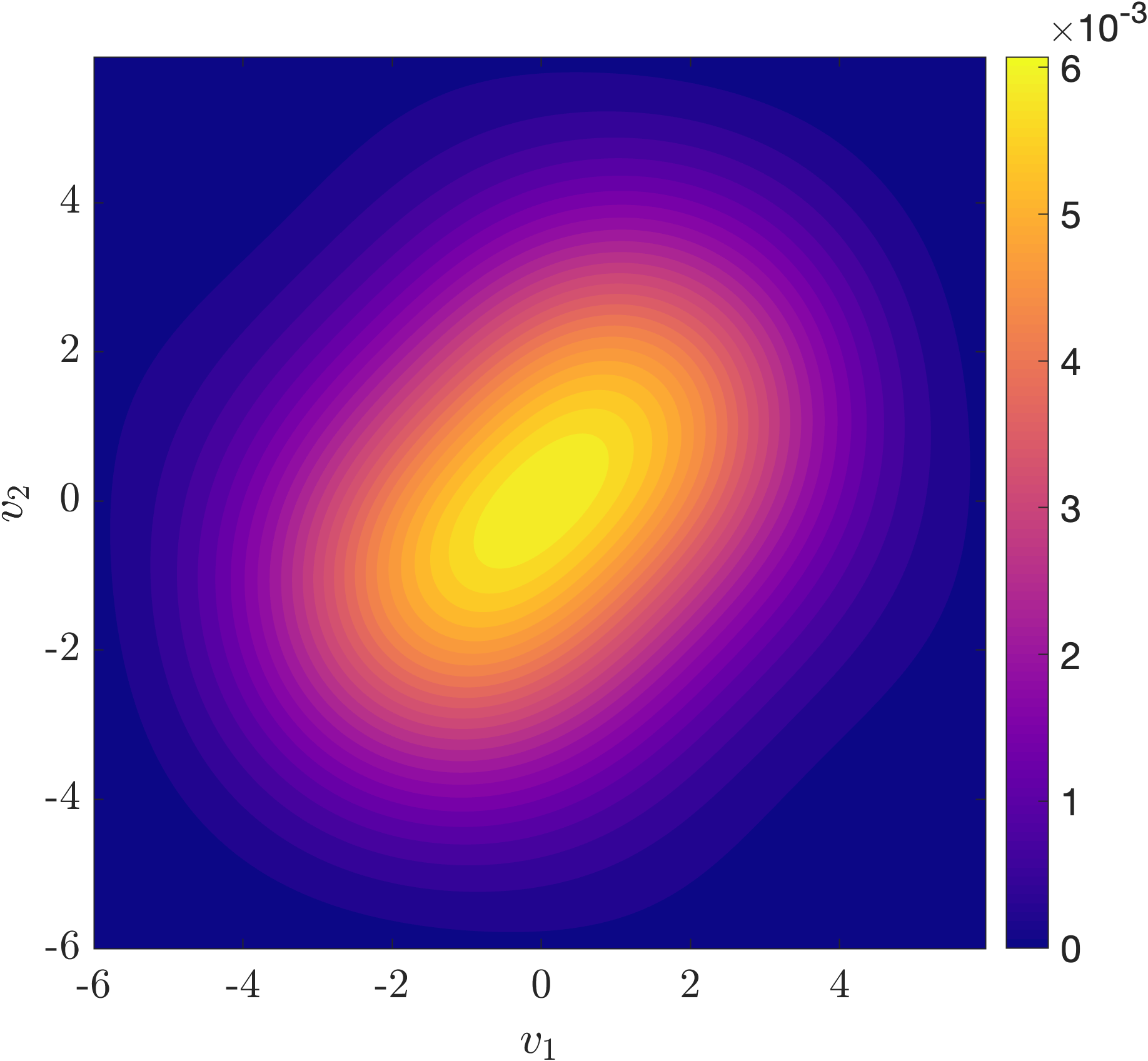}
        \caption{$t = 0.5$}
        \label{fig:dfp t = 0.5}
    \end{subfigure}
    \begin{subfigure}{0.32\linewidth}
        \centering
        \includegraphics[width=\linewidth]{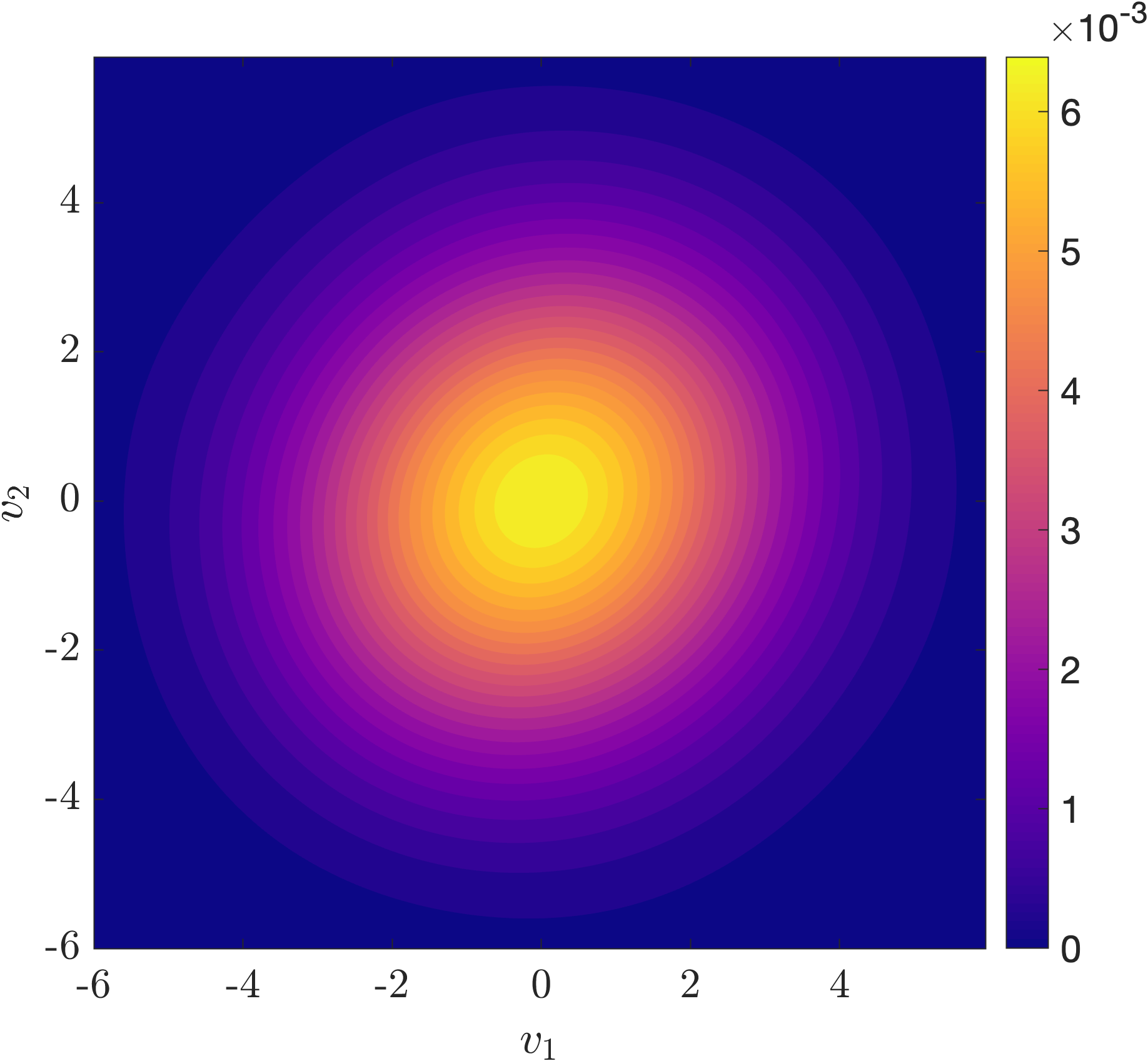}
        \caption{$t = 1.0$}
        \label{fig:dfp t = 1.0}
    \end{subfigure}
    
    \begin{subfigure}{0.32\linewidth}
        \centering
        \includegraphics[width=\linewidth]{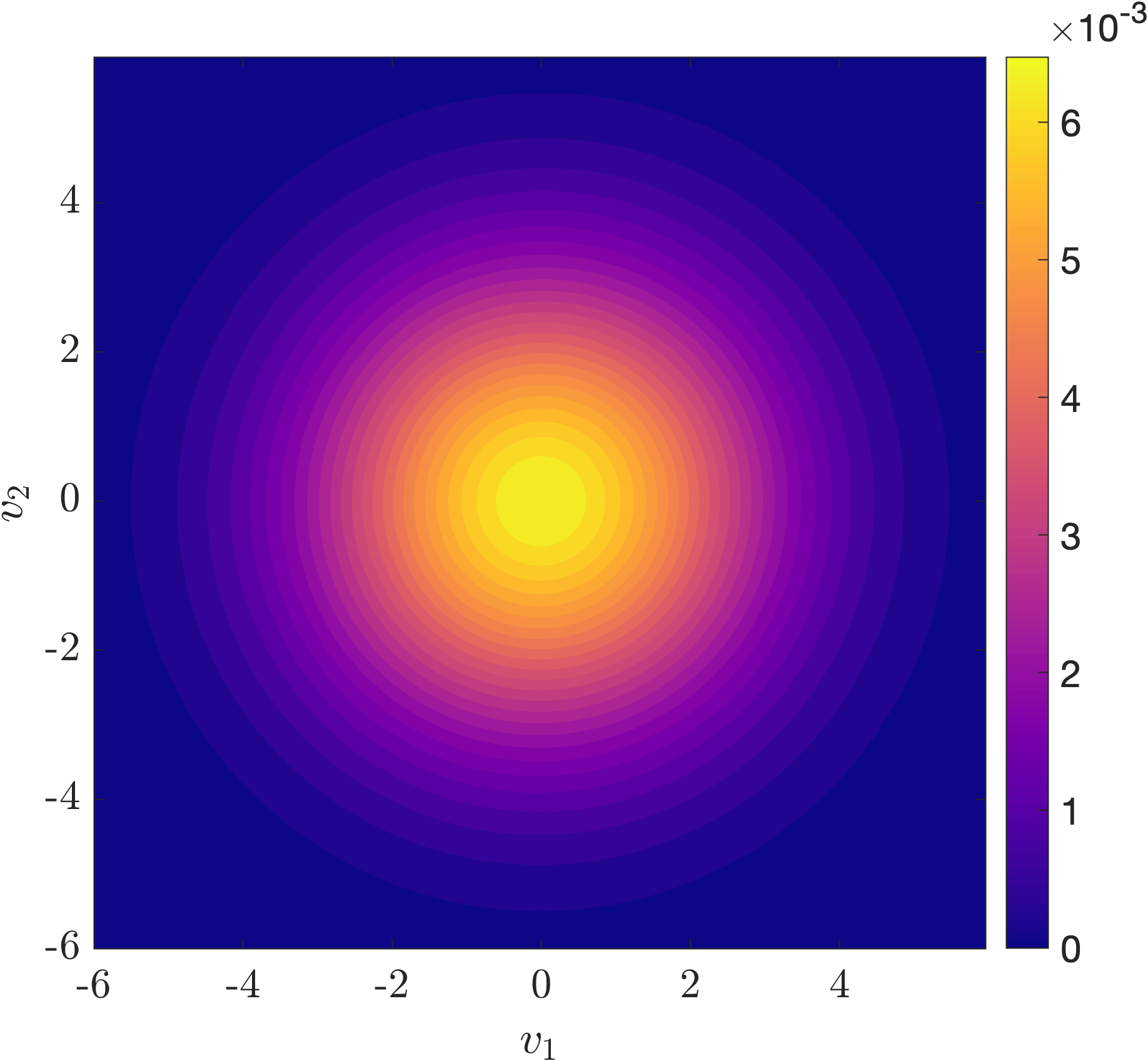}
        \caption{$t = 2.0$}
        \label{fig:dfp t = 2.0}
    \end{subfigure}
    \begin{subfigure}{0.32\linewidth}
        \centering
        \includegraphics[width=\linewidth]{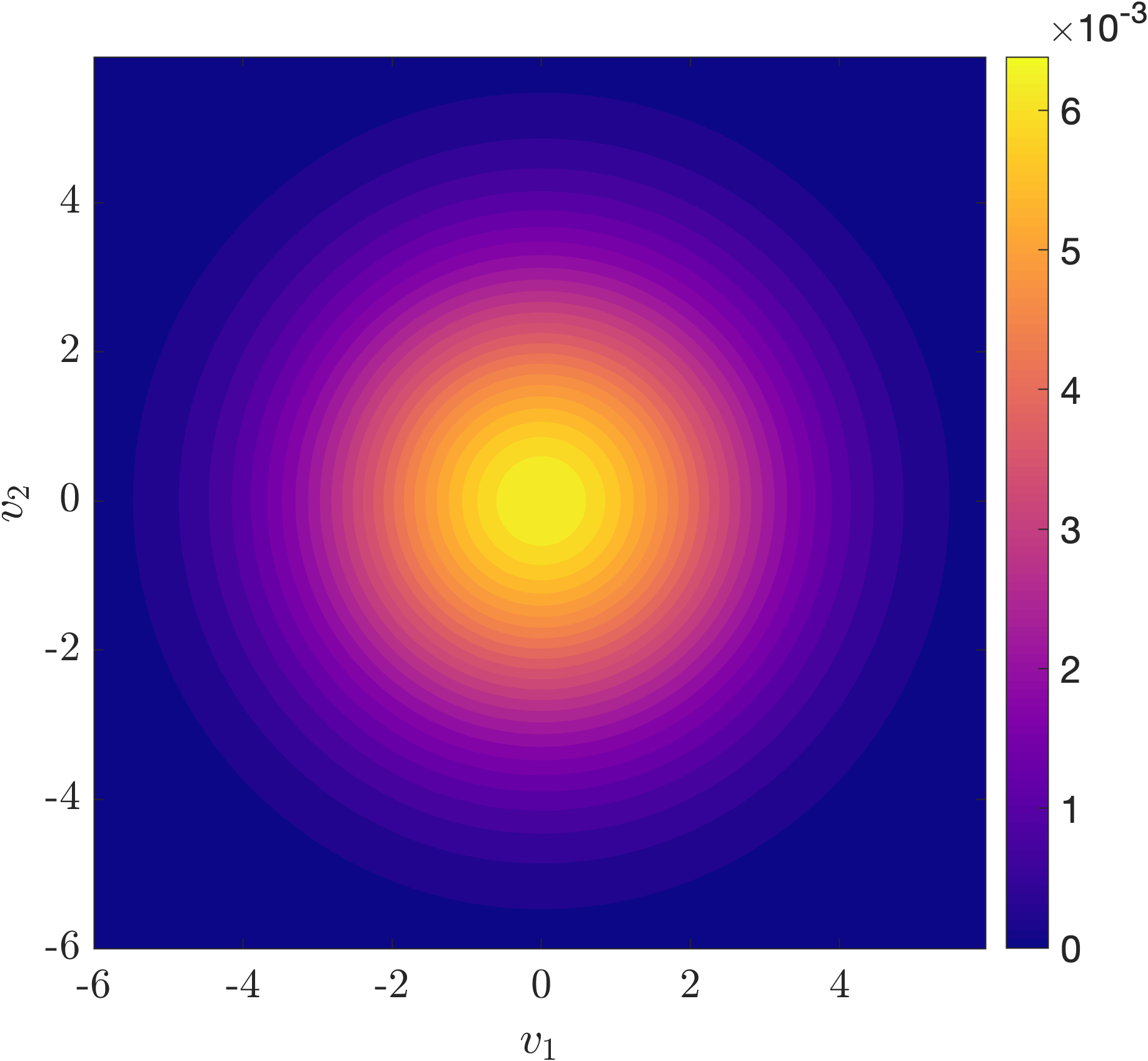}
        \caption{Equilibrium $f_{\infty}$}
        \label{fig:dfp equilibrium}
    \end{subfigure}
    \begin{subfigure}{0.32\linewidth}
        \centering
        \includegraphics[width=\linewidth]{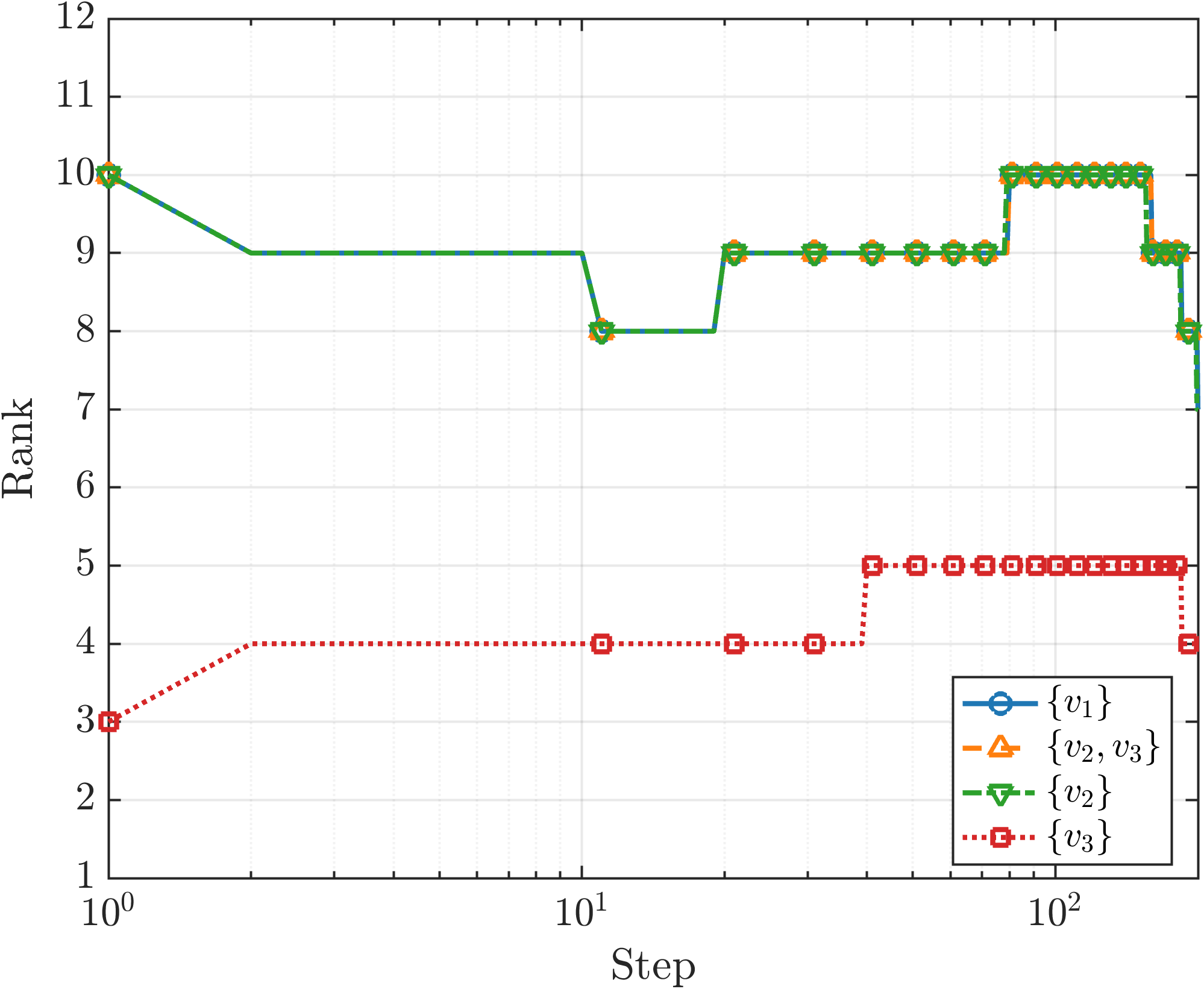}
        \caption{Rank vs step ($N = 16385$)}
        \label{fig:dfp rank data}
    \end{subfigure}

    \begin{subfigure}{0.325\linewidth}
        \centering
        \includegraphics[width=\linewidth]{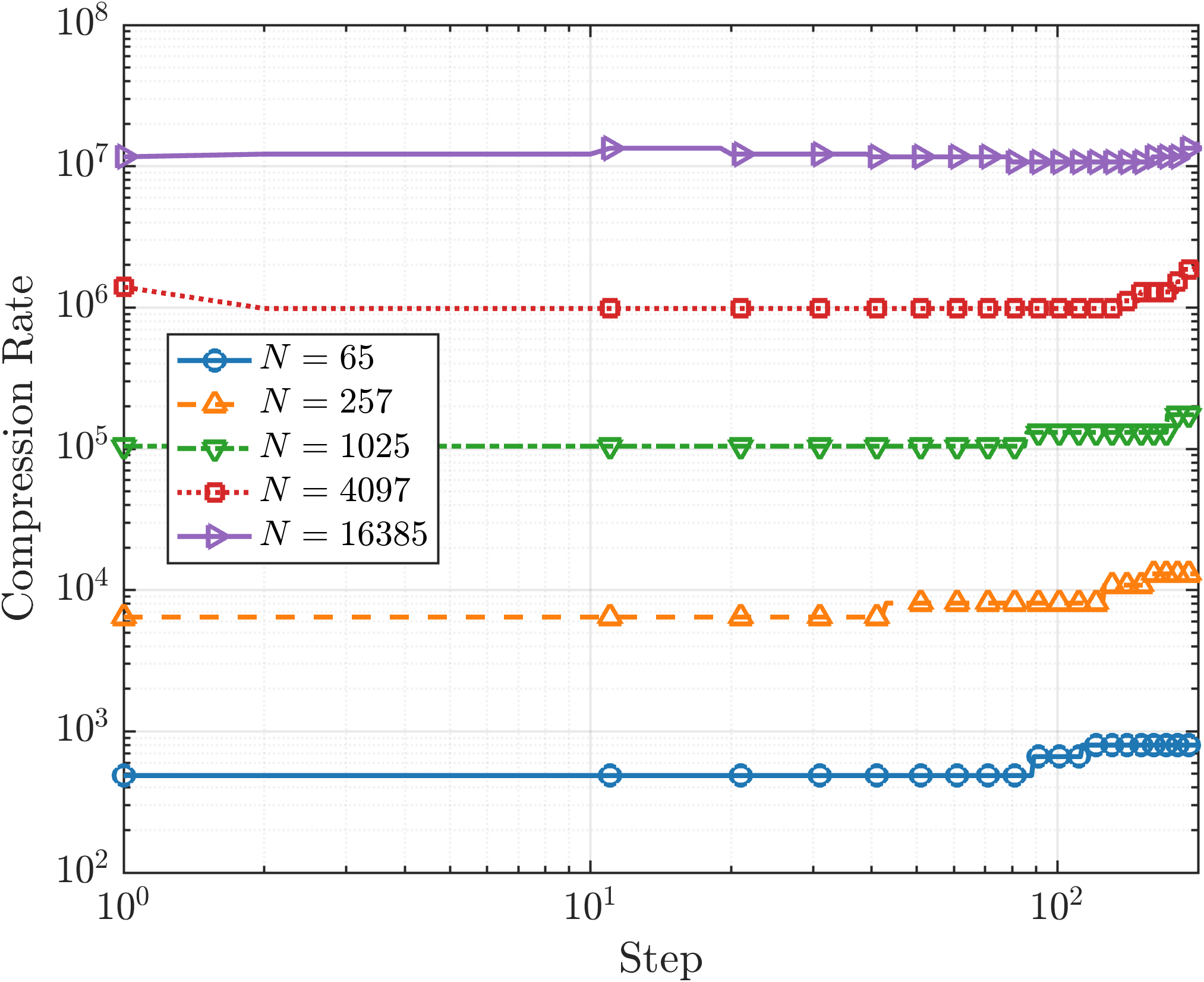}
        \caption{Compression rate vs step}
        \label{fig:dfp compression data}
    \end{subfigure}
    \begin{subfigure}{0.32\linewidth}
        \centering
        \includegraphics[width=\linewidth]{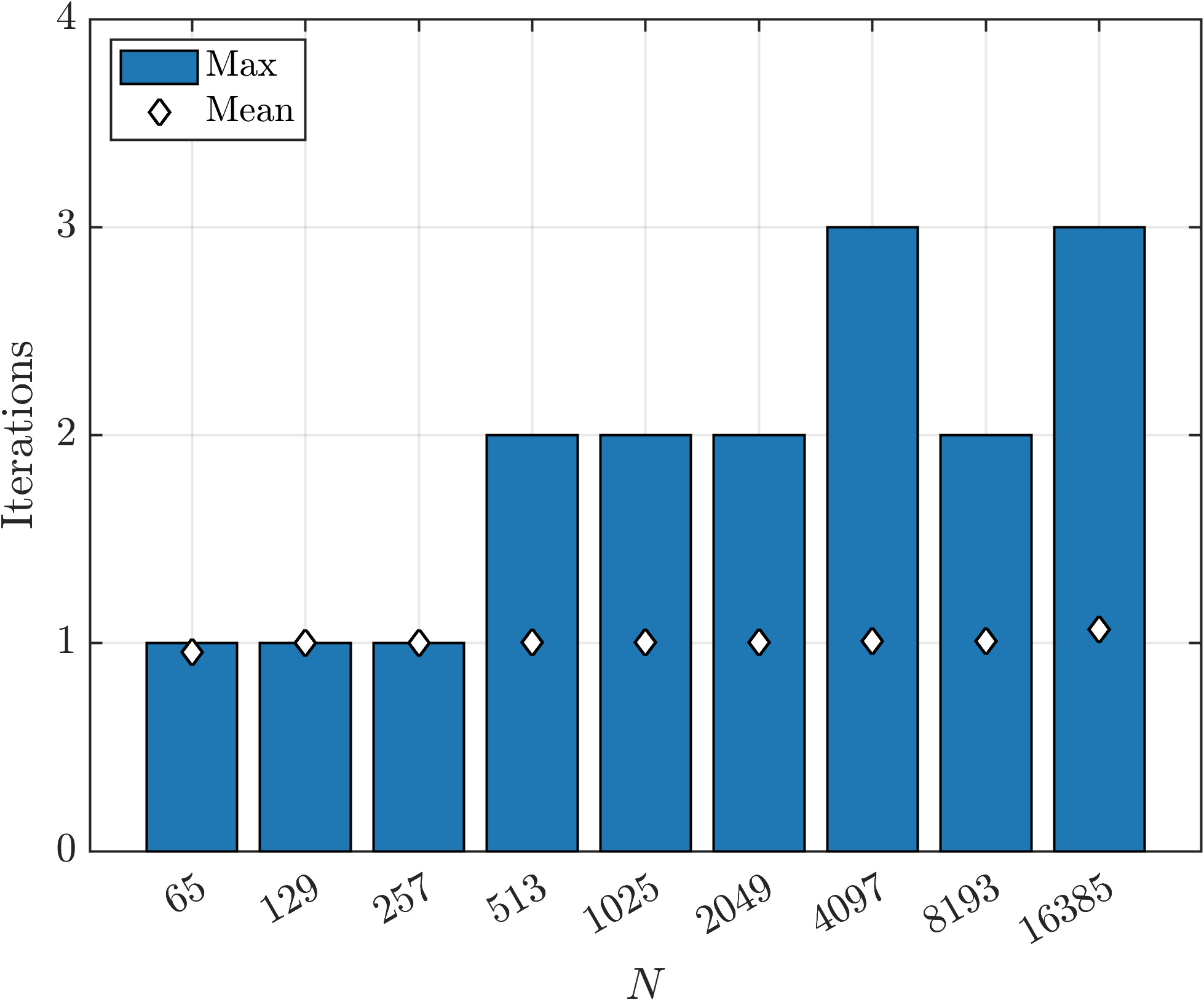}
        \caption{Iterations vs $N$}
        \label{fig:dfp mesh independence}
    \end{subfigure}
    \begin{subfigure}{0.3375\linewidth}
        \centering
        \includegraphics[width=\linewidth]{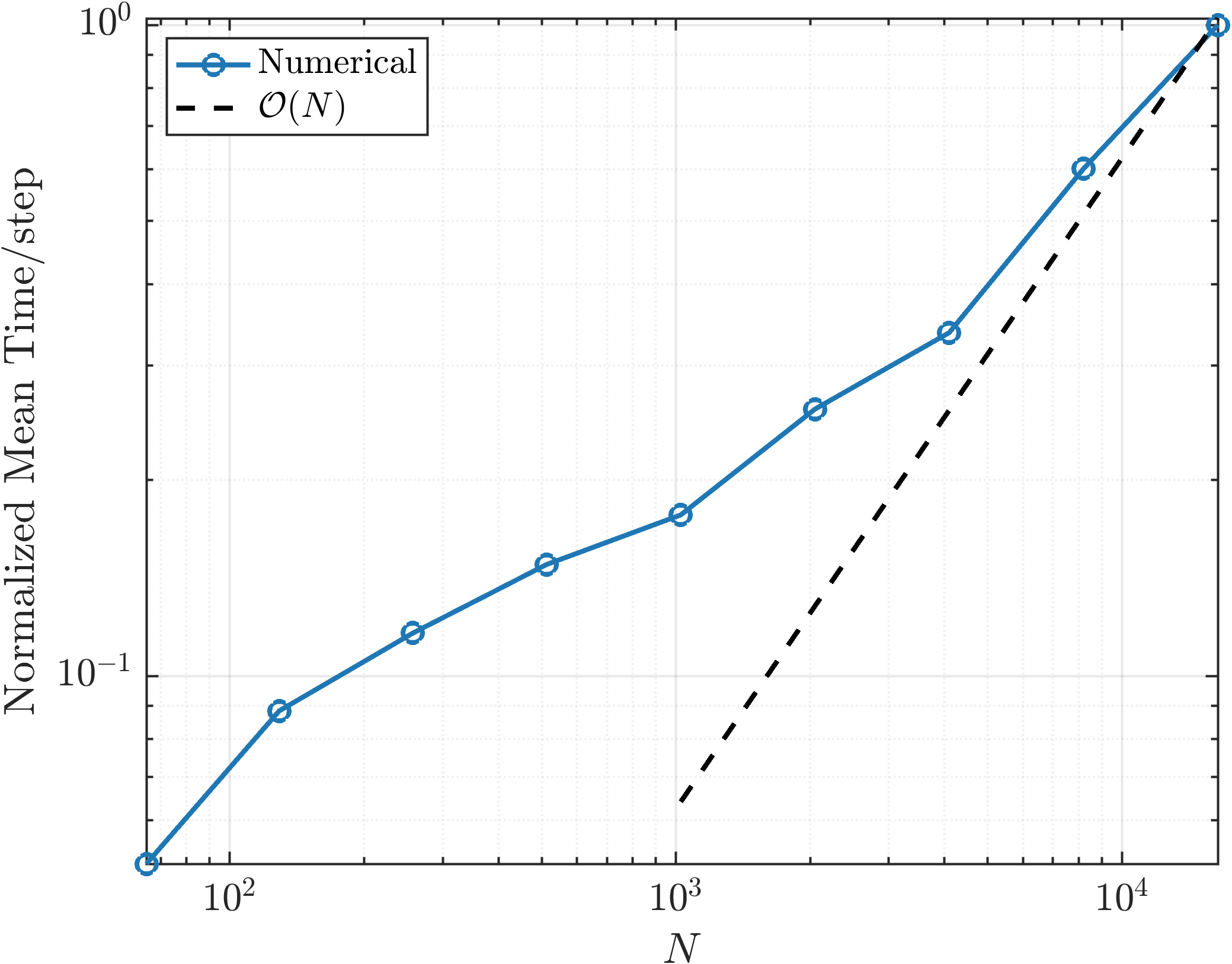}
        \caption{Complexity vs $N$}
        \label{fig:dfp complexity}
    \end{subfigure}
    \caption{(DFP problem)
    Simulation results obtained with right-preconditioned FGMRES using the GMG-V method. 
    (a -- d) Contour plots of the distribution function in the $v_1$--$v_2$ plane (with $v_3 = 0$) at selected times for a simulation with $N = 257$, using 25 equally spaced contours. The initial mixture of Maxwellians relaxes toward a single Maxwellian and is effectively at equilibrium by $t = 2$. 
    (e) The equilibrium distribution \eqref{eq:equilibrium distribution} is shown for reference.
    (f) Time history of the hierarchical ranks at each non-root node of the dimension tree for a fixed mode size $N = 16385$. The basis associated with $v_{3}$ has a lower rank than those for $v_{1}$ and $v_{2}$, reflecting the weaker dependence of the solution on $v_3$.
    (g) History of the compression rate, which increases slightly as the solution approaches equilibrium. 
    (h) Maximum and mean iteration counts for different values of the per-dimension mode size $N$, showing that the convergence depends only weakly on $N$.
    (i) Solver complexity as a function of the per-dimension mode size $N$. The mean time per step is normalized using the corresponding value at $N = 16385$.
    }
    \label{fig:DFP relaxation}
\end{figure}

We now consider the relaxation toward the equilibrium solution \eqref{eq:equilibrium distribution}. In this experiment we use a fixed time step size $\Delta t = 10^{-2}$ and set the final time to $T_f = 2.0$. For visualization, we use a velocity mesh with mode size $N = 257$ in each dimension. Figures \ref{fig:dfp t = 0} -- \ref{fig:dfp t = 2.0}  show contour plots of the distribution function in the $v_1$--$v_2$ plane (with $v_3 = 0$) at several representative times. The initial condition \eqref{eq:bimaxwellian initial data}, which consists of a mixture of Maxwellians, eventually coalesces into a single Maxwellian distribution. By $t = 2$, the solution has essentially reached its steady state and shows good qualitative agreement with the equilibrium distribution in Figure \ref{fig:dfp equilibrium}.

Using a fixed time step size $\Delta t = 10^{-2}$ and final time $T_{f} = 2.0$, we perform a sequence of mesh refinements in velocity space to study both the performance of the iterative solver and the resulting compression behavior. In each dimension, the mode size is chosen as $N = 2^{k} + 1$ with $k = 6, \ldots, 14$, ensuring that the grids remain nested across multigrid levels. The corresponding mesh spacing is $h = 12/(N-1)$. To quantify the stiffness of the problem, let $\lambda = T\Delta t/h^{2}$, which is a dimensionless parameter that characterizes the ratio between the implicit time scale and the parabolic stability restriction of an explicit scheme. For the time step and range of meshes considered here, $\lambda$ varies from $1.32$ to $8.70 \times 10^{2}$, indicating an increasingly stiff regime as the mesh is refined and highlighting the necessity of implicit time discretization.

We also examine the evolution of the hierarchical ranks for the case $N = 16385$, as shown in Figure~\ref{fig:dfp rank data}. The labels in the figure correspond to nodes in the dimension tree, while the rank of the root transfer tensor—associated with all dimensions—is omitted, as it is identically one by convention. The results indicate that the basis associated with the $v_{3}$ dimension consistently exhibits lower rank than those corresponding to $v_{1}$ and $v_{2}$, reflecting that the dominant dynamics occur in the $v_{1}$–$v_{2}$ plane. As the solution approaches equilibrium, the hierarchical ranks decrease appreciably. This behavior is consistent with the concurrent increase in compression observed in Figure~\ref{fig:dfp compression data}, which examines the compression rate achieved by the low-rank representation during the simulation. The compression rate increases with $N$ because the storage for the full-grid scheme grows as $N^3$, whereas the storage for the adaptive-rank scheme grows linearly with respect to $N$. Furthermore, we observe a slight increase in each case as the solution approaches equilibrium.

In Figures \ref{fig:dfp mesh independence} and \ref{fig:dfp complexity}, we illustrate the scalability of the solver with respect to the per-dimension mode size $N$. In Figure \ref{fig:dfp mesh independence}, we report the mean and maximum number of FGMRES iterations per time step over the full simulation for each $N$. The iteration counts remain essentially unchanged as the mesh is refined, indicating that the convergence behavior is effectively independent of $N$. Figure \ref{fig:dfp complexity} shows that the normalized average time per step scales approximately as $\mathcal{O}(N)$, yielding near-linear overall complexity.

\subsection{Allen--Cahn Equation}
\label{subsec:allen_cahn}


The next example we consider is the 3D Allen--Cahn equation on the unit cube $\Omega_{\mathbf{x}} = [0,1]^{3}$
\begin{equation}
\label{eq:Allen-Cahn system}
\left\{
\begin{aligned}
    &\partial_{t}u = \epsilon^{2}\Delta_{\mathbf{x}} u - f(u),
    &&(\mathbf{x},t) \in \Omega_{\mathbf{x}} \times (0,T_{f}], \\
    &u(\mathbf{x},0) = u_{0}(\mathbf{x}),
    &&\mathbf{x} \in \Omega_{\mathbf{x}}.
\end{aligned}
\right.
\end{equation}
subject to periodic boundary conditions. 

The Allen--Cahn equation is a nonlinear reaction–diffusion equation arising in models of phase separation in binary mixtures \cite{allen1979microscopic}. The parameter $0 < \epsilon$ represents the characteristic width of the interface separating phases, while the nonlinear term
\begin{equation*}
    f(u) = u^{3} - u,
\end{equation*}
corresponds to a double-well potential with stable states at $u = \pm 1$. The isosurface $u = 0$ is naturally interpreted as the interface between phases.

Equation \eqref{eq:Allen-Cahn system} is the $L^{2}$-gradient flow of the Ginzburg–Landau energy functional
\begin{equation}
    \label{eq:Ginzburg-Landau energy}
    \mathcal{E}[u] = \int_{\Omega_{\mathbf{x}}} \left( \frac{\epsilon^{2}}{2}|\nabla u|^{2} + \frac{1}{4}(u^{2} - 1)^{2} \right)\,d\mathbf{x},
\end{equation}
which satisfies the dissipation property
\begin{equation}
    \label{eq:Allen-Cahn energy dissipation}
    \frac{d\mathcal{E}[u]}{dt} \leq 0.
\end{equation}
As a result, the dynamics drive the system toward configurations that reduce interfacial curvature and minimize the total energy.

In this experiment, we investigate whether low-rank methods can capture these dynamics while maintaining compression. The initial data in \eqref{eq:Allen-Cahn system} is taken as a low-rank perturbation of the unstable state $u=0$, given by
\begin{equation}
    \label{eq:Allen-Cahn spinodal IC}
    u_{0}(\mathbf{x}) = \eta\left( \prod_{\mu=1}^{3}\cos \left(k_{\mu}x_{\mu}\right) + \sin \left(2k_{1}x_{1} + k_{1}x_{2}\right) + \cos \left(k_{1}x_{2} + k_{1}x_{3}\right)  \right),
\end{equation}
with parameters
\begin{equation*}
    \eta = 0.3, \quad k_{1} = 2\pi, \quad k_{2} = -4\pi, \quad k_{3} = 6\pi.
\end{equation*}

Discretizing the PDE \eqref{eq:Allen-Cahn system} in time using the backward Euler method gives the update
\begin{equation*}
    u^{n+1} - u^{n} + \Delta t \Big[ -\epsilon^{2} \Delta_{\mathbf{x}}u^{n+1} + f \left( u^{n+1} \right) \Big] = 0.
\end{equation*}
The Laplacian is discretized with mode-wise periodic finite difference approximations with second-order accuracy. This leads to the nonlinear system
\begin{equation}
    \label{eq:Allen-Cahn nonlinear system}
    \mathcal{F}(u)
    :=
    u - u^{n}
    + \Delta t \left[
    -\epsilon^{2}
    \sum_{\mu=1}^{3} D_{x_{\mu}x_{\mu}} u
    + f(u)
    \right] = 0,
\end{equation}
Linearizing this system about $u=u^{*}$ gives the Jacobian-vector product
\begin{equation*}
    \mathcal{J}(u^{*})v
    =
    v
    +
    \Delta t \left[
    -\epsilon^{2}
    \sum_{\mu=1}^{3} D_{x_{\mu}x_{\mu}} v
    +
    f'(u^{*}) \odot v
    \right].
\end{equation*}
Here
\begin{equation*}
    f'(u^{*}) = 3 u^{*} \odot u^{*} - \bigotimes_{\mu=1}^{d} \mathbf{1}_{\mu},
\end{equation*}
where $\mathbf{1}_{\mu}$ denotes the vector of length $n_{\mu}$ whose entries are all equal to one. Alternatively, one can compute these Jacobian vector products using the finite differencing in \eqref{eq:Jacobian-vector product approximation}. In this experiment, we use a fixed time step size $\Delta t = 10^{-2}$ and set the final time to $T_f = 2.0$.

The nonlinear system \eqref{eq:Allen-Cahn nonlinear system} in each time step is solved using the inexact Newton method in Section \ref{subsec:inexact Newton tensor methods}. The nonlinear iteration is terminated when the residual satisfies either the absolute or relative stopping criterion \eqref{eq:Newton stopping residual tolerances}, with $\tau_{\mathrm{abs}} = 10^{-12}$ and $\tau_{\mathrm{rel}} = 10^{-4}$. The methods also terminate if the Newton correction tolerances satisfy \eqref{eq:Newton stopping update tolerances}, with $\xi_{\mathrm{abs}} = 10^{-12}$ and $\xi_{\mathrm{rel}} = 10^{-10}$. For the inner linear iterations, we use the adaptive relative tolerance defined by \eqref{eq:Inner Newton forced relative tolerance} and \eqref{eq:forcing term}, with $\gamma_{\min} = 10^{-4}$ and $\gamma_{\max} = 1/2$. Since the nonlinear solve is sensitive to truncation errors, we use the more conservative truncation tolerances $\epsilon_{\mathrm{abs}} = \epsilon_{\mathrm{rel}} = 10^{-6}$. The algorithm uses the Armijo backtracking \eqref{eq:Armijo backtracking} with $c = 10^{-4}$ and a scalar reduction factor of $0.5$. The inner linear systems are solved using FGMRES with adaptive tolerances set according to \eqref{eq:Inner Newton forced relative tolerance} with GMG-V preconditioning based on the same tolerances. The preconditioner approximately inverts the operator
\begin{equation*}
    \widetilde{\mathcal{J}}v = v - \epsilon^{2}\Delta t \sum_{\mu=1}^{3} D_{x_{\mu}x_{\mu}} v,
\end{equation*}
using two V-cycles and uses 10 damped Jacobi iterations with $\omega=0.7$ as a smoother. On level $\ell$, the elementwise inverse of the diagonal operator is a scalar multiple of the identity:
\begin{equation*}
    \left(\mathcal{D}^{(\ell)}\right)^{-1} = \left[ 1 + 2\epsilon^{2}\Delta t
     \sum_{\mu=1}^d \frac{1}{\left(h_{\mu}^{(\ell)}\right)^2} \right]^{-1} \mathcal{I}.
\end{equation*}
An advantage of this preconditioner is that it does not depend on the Newton state, so it can be constructed once and reused during the time integration.

\begin{figure}[!htbp]
    \centering
    \begin{subfigure}{0.315\linewidth}
        \centering
        \includegraphics[width=\linewidth]{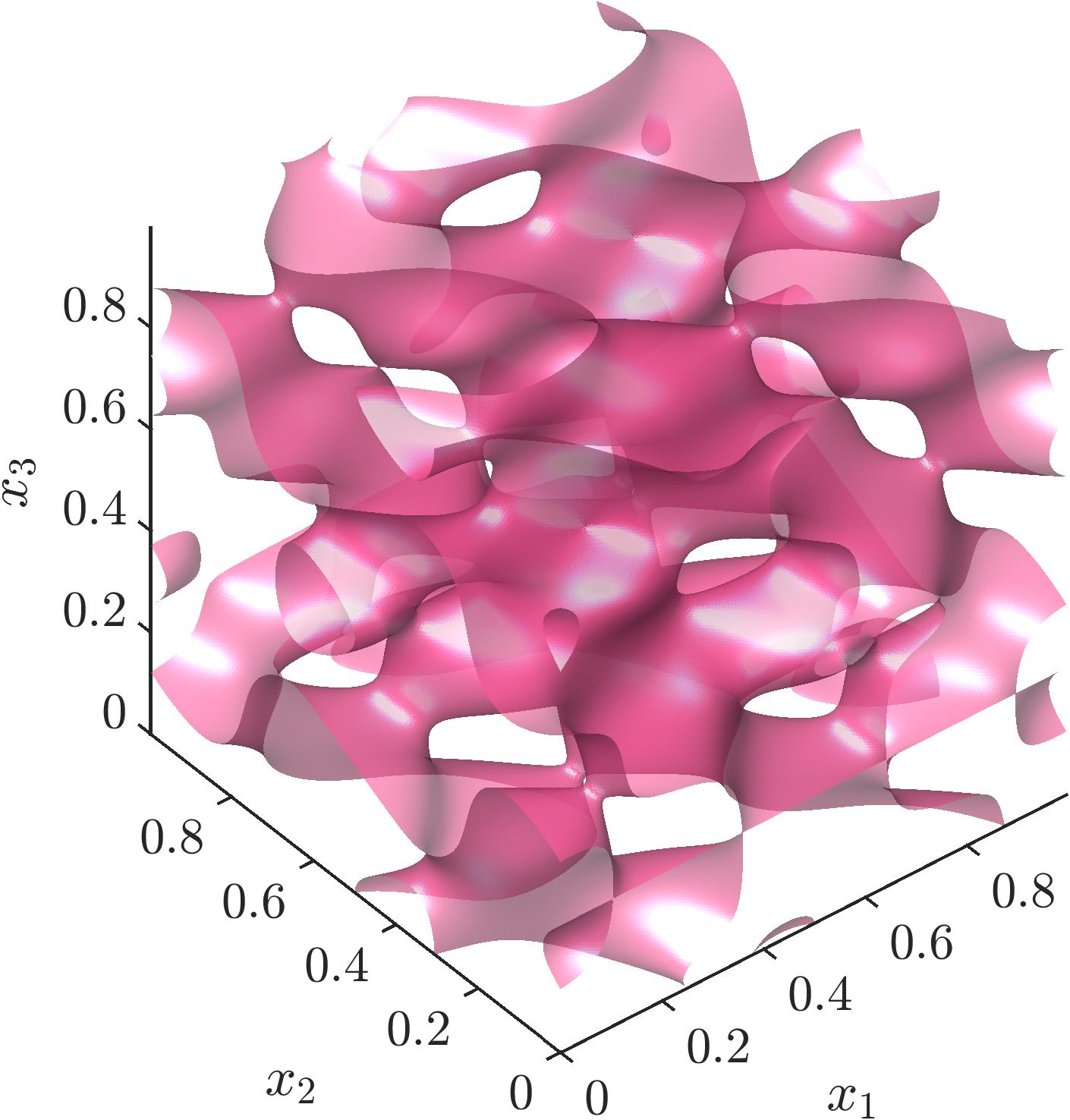}
        \caption{$t = 0$}
        \label{fig:ac t = 0}
    \end{subfigure}
    \begin{subfigure}{0.315\linewidth}
        \centering
        \includegraphics[width=\linewidth]{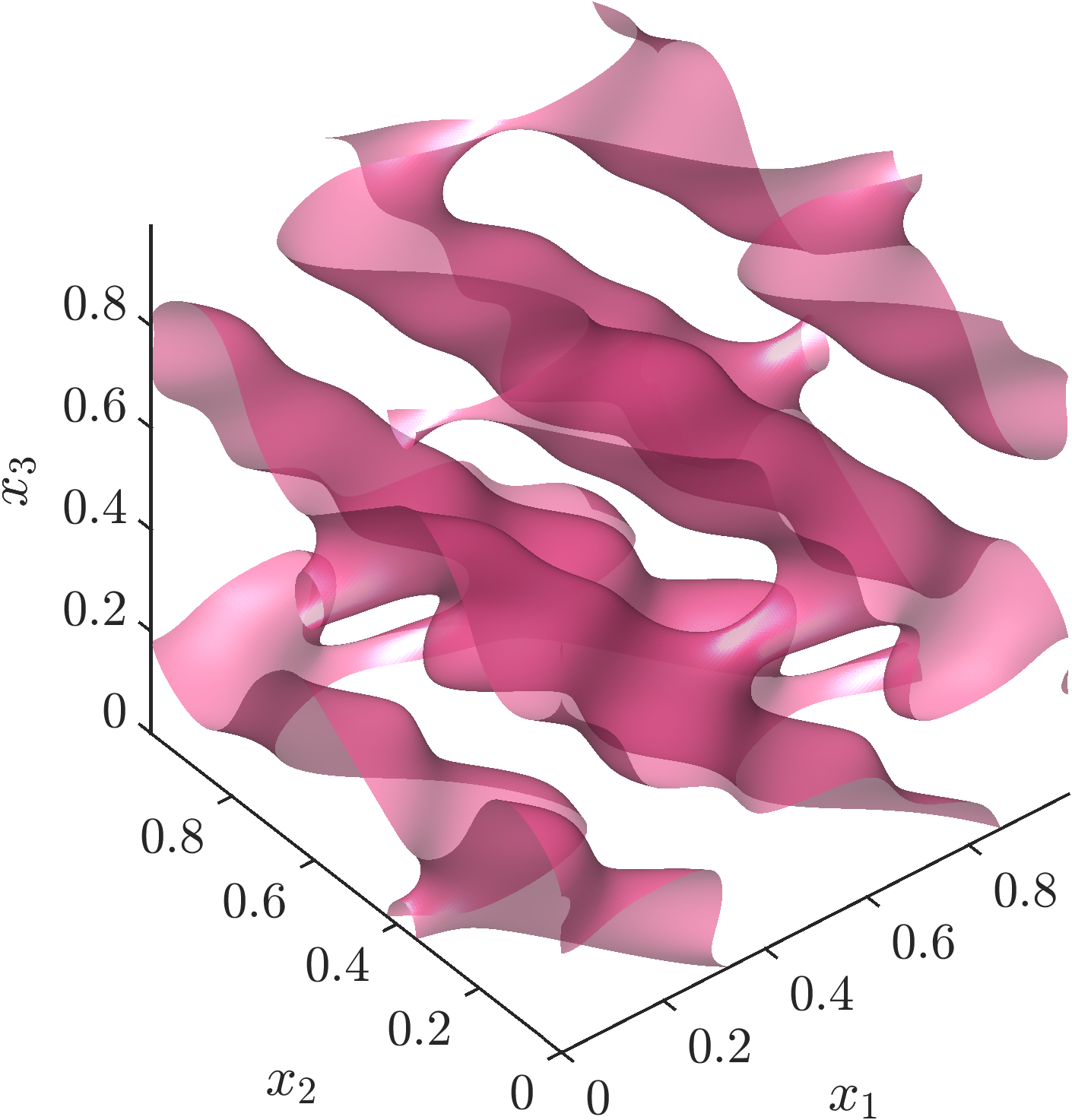}
        \caption{$t = 0.2$}
        \label{fig:ac t = 0.2}
    \end{subfigure}
    \begin{subfigure}{0.32\linewidth}
        \centering
        \includegraphics[width=\linewidth]{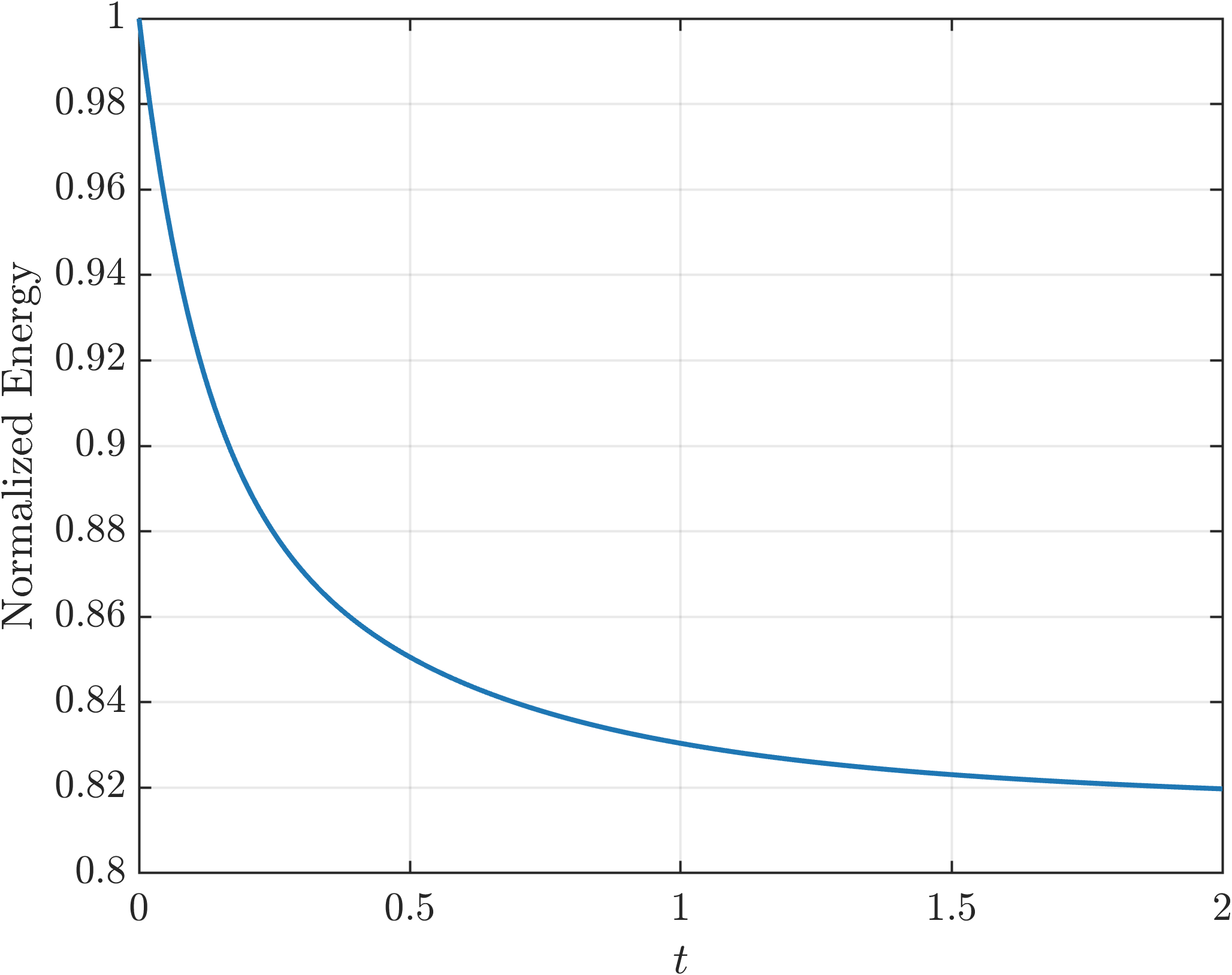}
        \caption{Normalized energy ($N = 257$)}
        \label{fig:ac energy}
    \end{subfigure}

    \begin{subfigure}{0.315\linewidth}
        \centering
        \includegraphics[width=\linewidth]{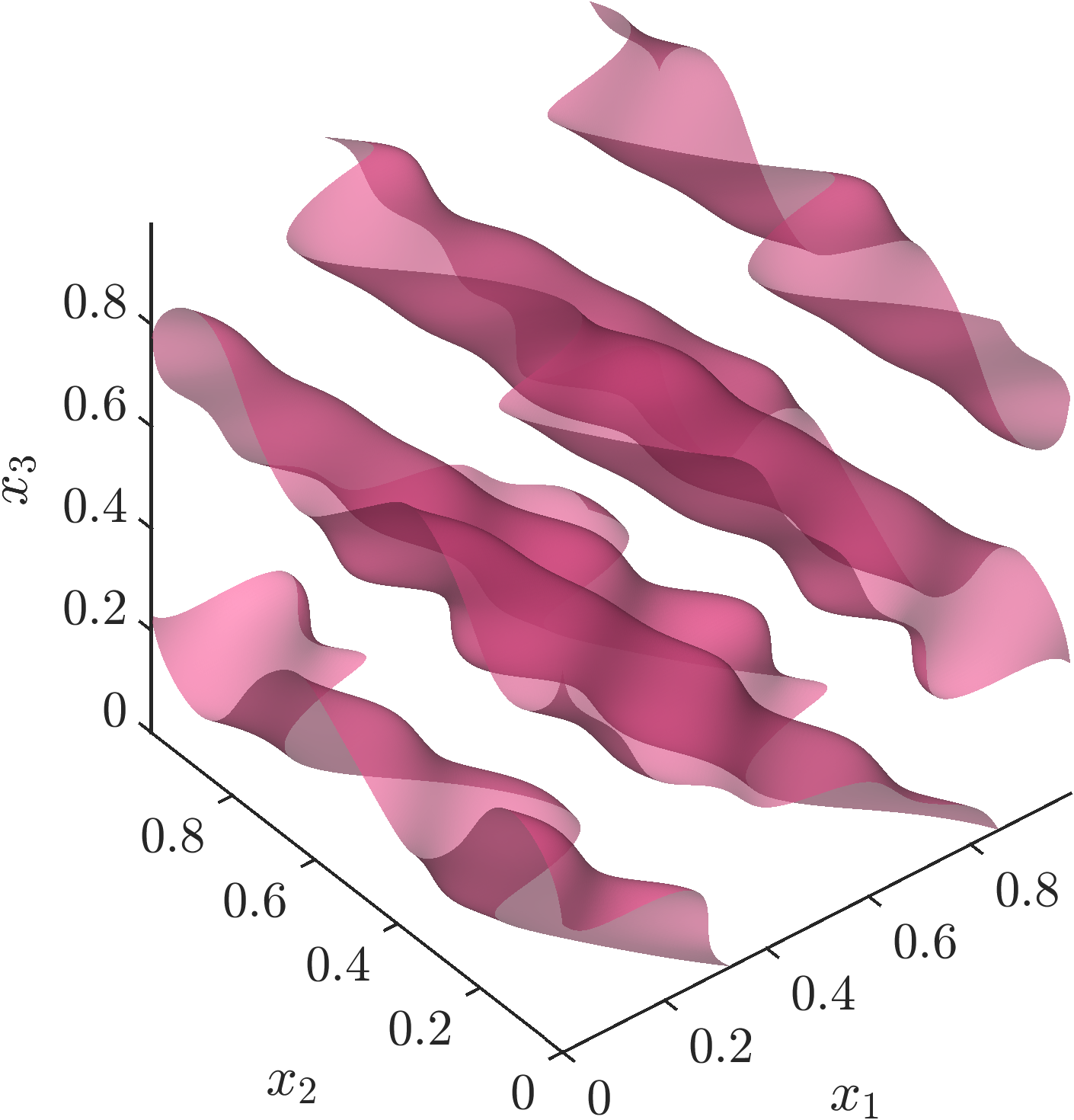}
        \caption{$t = 0.3$}
        \label{fig:ac t = 0.3}
    \end{subfigure}
    \begin{subfigure}{0.315\linewidth}
        \centering
        \includegraphics[width=\linewidth]{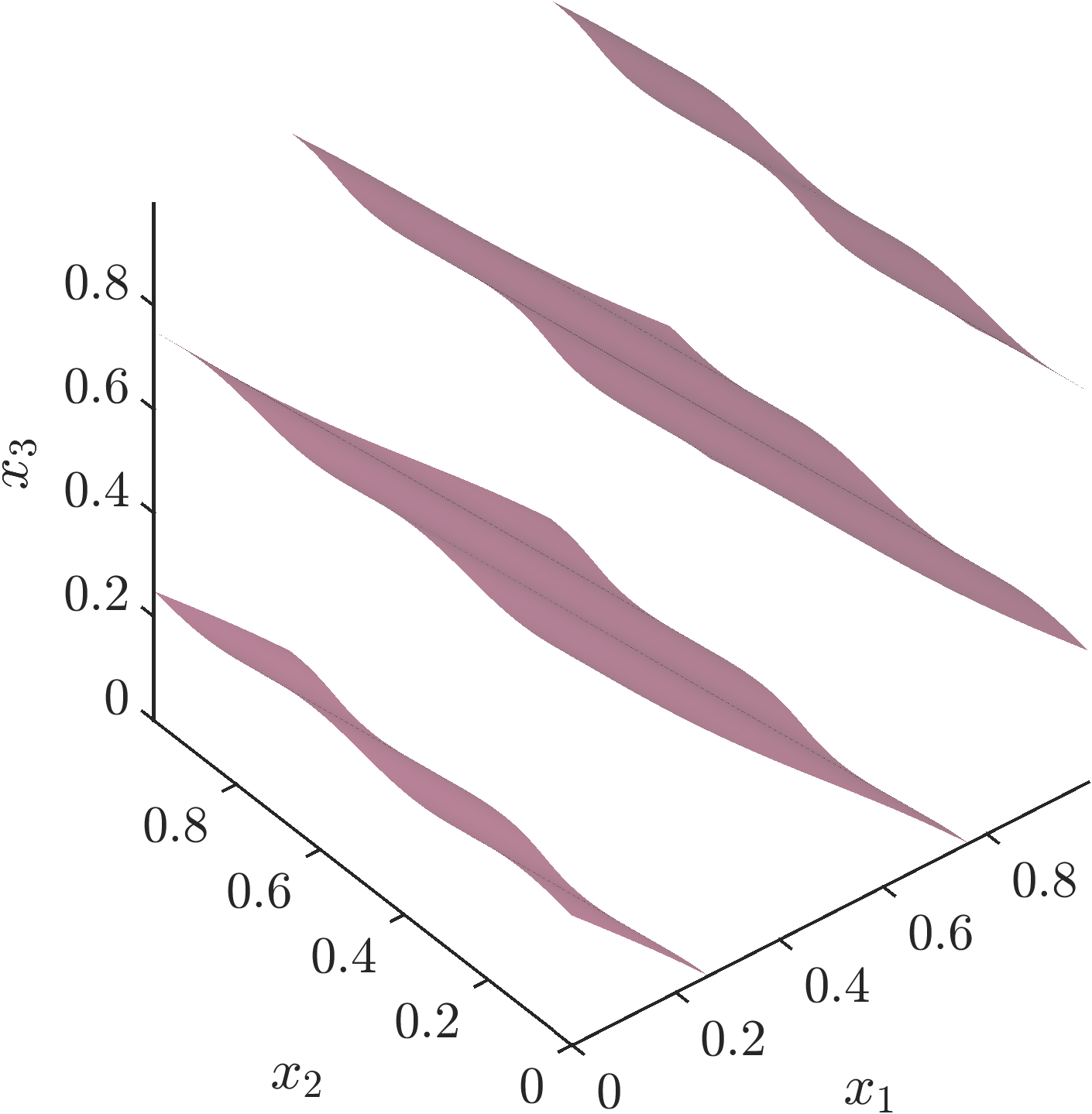}
        \caption{$t = 2.0$}
        \label{fig:ac t = 2.0}
    \end{subfigure}
    \begin{subfigure}{0.32\linewidth}
        \centering
        \includegraphics[width=\linewidth]{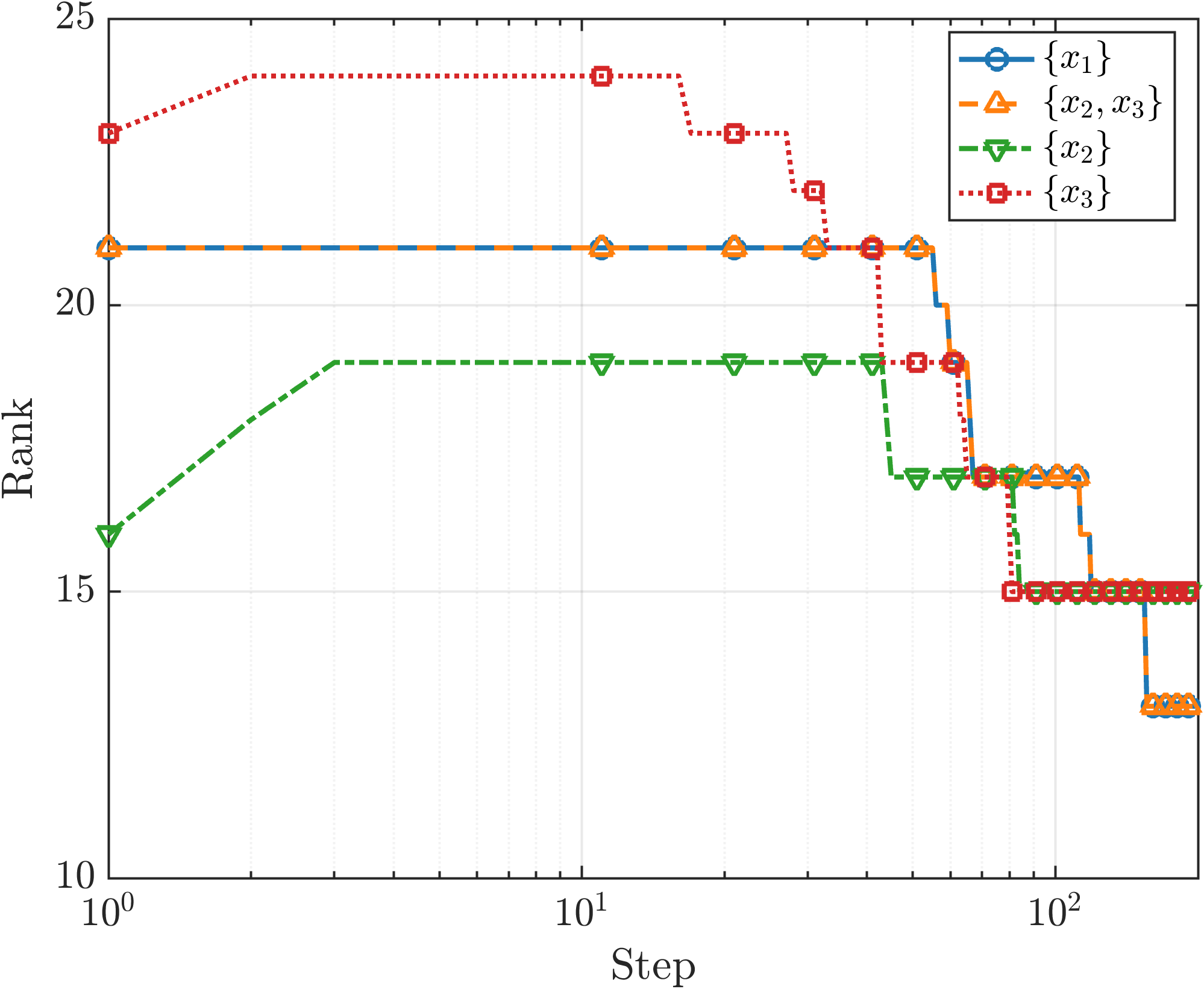}
        \caption{Rank vs step ($N = 257$)}
        \label{fig:ac rank}
    \end{subfigure}

    \begin{subfigure}{0.3285\linewidth}
        \centering
        \includegraphics[width=\linewidth]{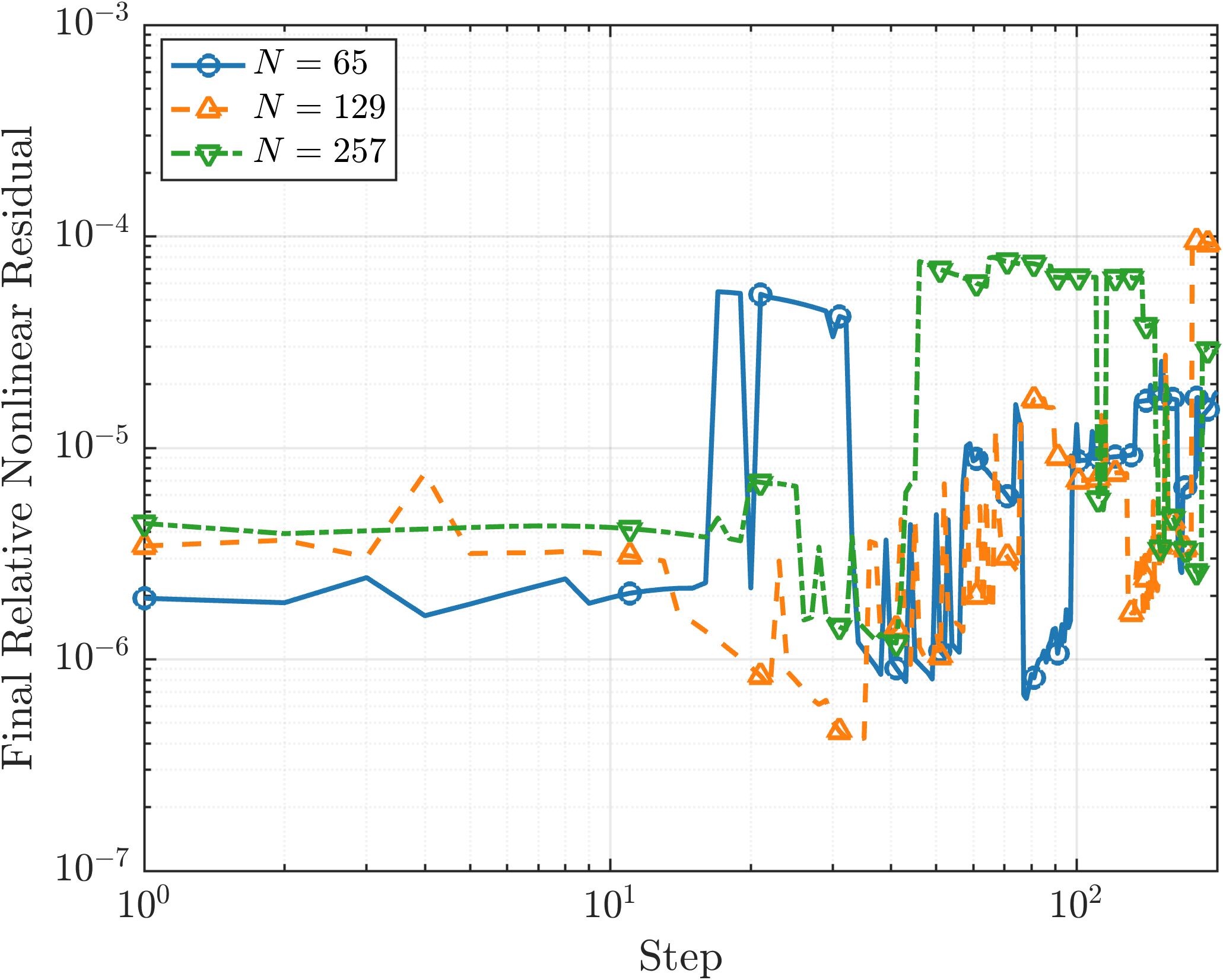}
        \caption{Final relative residual vs step}
        \label{fig:ac nonlinear residual}
    \end{subfigure}
    \begin{subfigure}{0.325\linewidth}
        \centering
        \includegraphics[width=\linewidth]{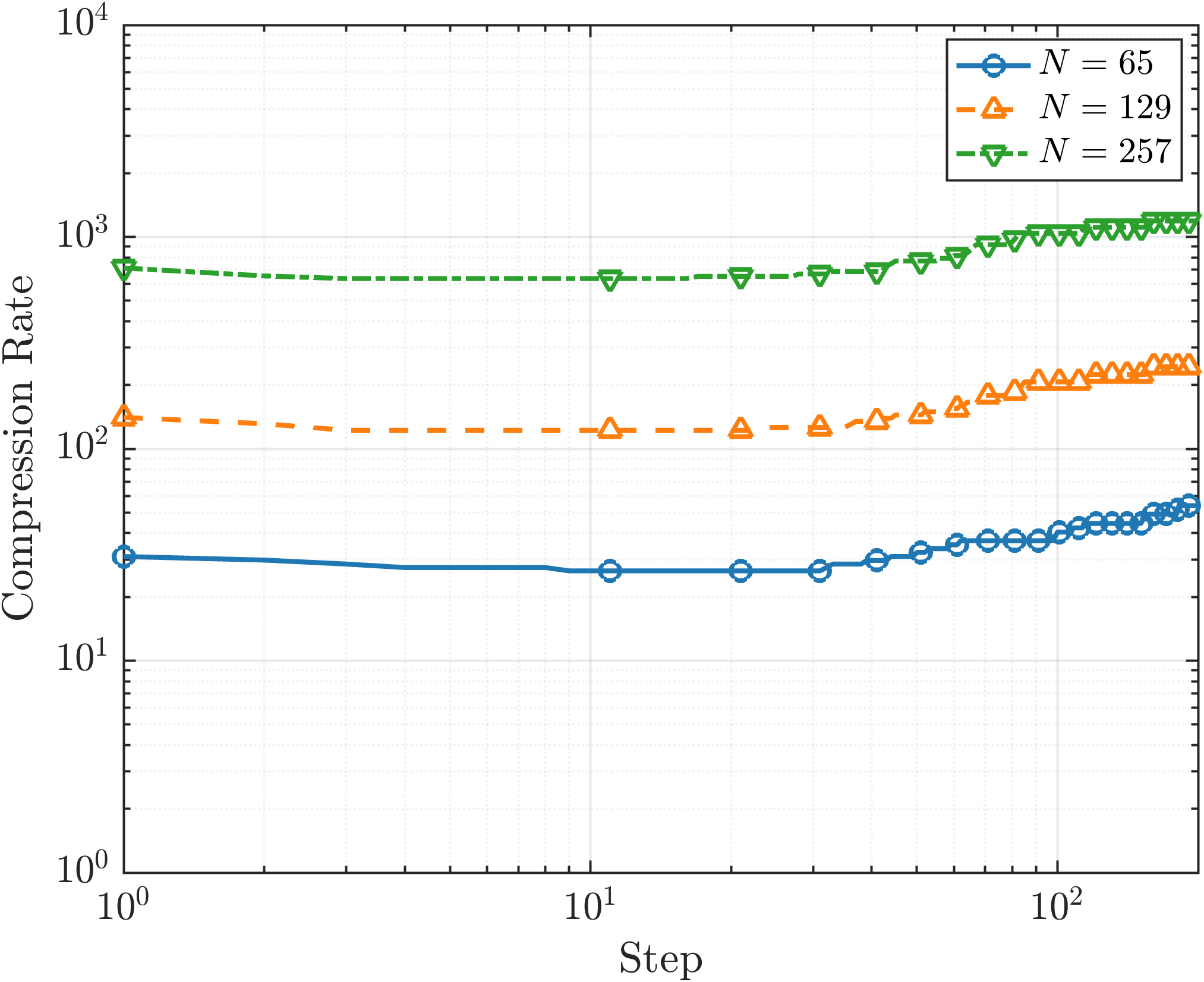}
        \caption{Compression rate vs step}
        \label{fig:ac compression}
    \end{subfigure}
    \begin{subfigure}{0.3175\linewidth}
        \centering
        \includegraphics[width=\linewidth]{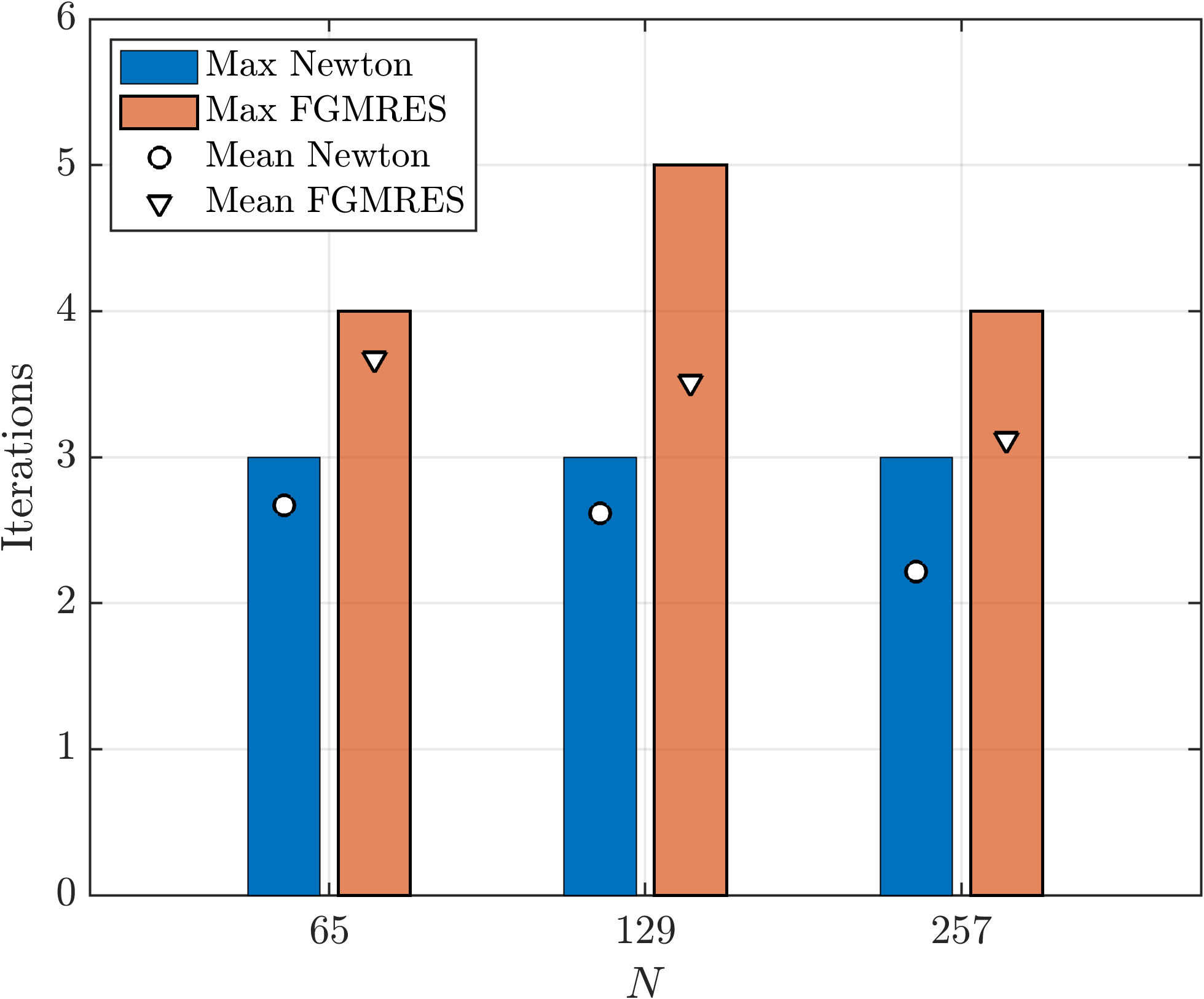}
        \caption{Iterations vs $N$}
        \label{fig:ac iterations}
    \end{subfigure}
    \caption{(Allen--Cahn problem) Simulation results obtained using the inexact Newton method with right-preconditioned FGMRES. (a), (b), (d), and (e) Snapshots of the isosurface \(u=0\) for \(N=257\), illustrating the coarsening dynamics and curvature driven flow toward the steady state. (c) The monotone decay of the normalized energy, which is consistent with equation \eqref{eq:Allen-Cahn energy dissipation}. (f) The hierarchical ranks for each non-root node in the dimension tree. The larger rank associated with \(x_{3}\) reflects the higher-frequency content specified by \(k_{3}\) in \eqref{eq:Allen-Cahn spinodal IC}. (g) and (h) The final relative nonlinear residuals and compression rates for several values of \(N\), respectively. In all cases, the nonlinear solves converge within the prescribed tolerance, while the compression rate increases as the solution approaches equilibrium. (i) Summary of the Newton and FGMRES iteration counts across mesh sizes, with bars denoting maxima over the simulation and markers denoting averages. The iteration counts remain largely unchanged as \(N\) is varied.
    }
\end{figure}

Figures \ref{fig:ac t = 0}, \ref{fig:ac t = 0.2}, \ref{fig:ac t = 0.3}, and \ref{fig:ac t = 2.0} show the evolution of the interface $u=0$ at different stages of the simulation with $N = 257$. The initially complex structure of the interface rapidly coalesces into smooth, sheet-like structures that are not aligned with the coordinate axes. By $t=2$, the solution has effectively reached a steady state, which is reflected by the energy decay shown in Figure \ref{fig:ac energy}. The energy, normalized with respect to the initial condition, decays monotonically, which is consistent with equation \eqref{eq:Allen-Cahn energy dissipation}. In our implementation, we compute the discrete analogue of \eqref{eq:Ginzburg-Landau energy} as
\begin{equation*}
    \mathcal{E}[u]
    \approx
    \left(\prod_{\mu=1}^{3} h_{\mu} \right)
    \left[
    \frac{\epsilon^2}{2}
    \sum_{\mu=1}^{d} \|D_{x_\mu}u\|_{F}^2
    +
    \frac{1}{4}
    \Big\| u\odot u - \bigotimes_{\mu=1}^{d} \mathbf{1}_{\mu} \Big\|_{F}^{2}
    \right],
\end{equation*}
where the spatial derivatives are approximated by periodic central differences. Figure \ref{fig:ac rank} presents an alternative form of this information in the form of hierarchical ranks. At early steps, the rank of the solution increases to accommodate the topology changes in the solution. However, as the solution approaches steady state, the ranks show a uniform decrease as the structure of the solution simplifies. We observe that the hierarchical rank associated with the leaf node corresponding to $x_{3}$ is notably higher than those associated with $x_{1}$ and $x_{2}$. This reflects the higher-frequency content associated with $k_3$ in the initial condition \eqref{eq:Allen-Cahn spinodal IC}.

Figures \ref{fig:ac nonlinear residual} -- \ref{fig:ac iterations} present several simulation diagnostics obtained with different per-dimension mode sizes $N$. Figure \ref{fig:ac nonlinear residual} shows the relative nonlinear residual obtained from the final nonlinear iterate at each time step. In all cases, the residuals remain within the prescribed tolerance, indicating reliable convergence of the nonlinear solver. Figure \ref{fig:ac compression} shows the corresponding compression histories. Compared with the previous test problems, the compression achieved in this experiment is more modest. This is due to the more complex solution structure and the fact that the evolving interfaces are not aligned with the coordinate axes of the computational domain. Nevertheless, the proposed method still achieves substantial compression in several cases, even during periods in which the topology of the solution changes rapidly. As the solution approaches steady state, the compression rate increases, reflecting the emergence of simpler solution structures. Finally, Figure~\ref{fig:ac iterations} reports iteration counts for the inexact Newton solver and the right-preconditioned FGMRES method. The number of Newton iterations remains at most three for each value of \(N\), including during the early stages of the simulation when the solution undergoes significant topological changes. The FGMRES iteration counts also remain small, demonstrating the effectiveness of the multilevel preconditioner.

\section{Conclusion}
\label{sec:conclusion}


In this work, we developed a unified framework for solving high-dimensional PDEs using adaptive low-rank tensor representations in combination with multilevel preconditioning. The approach uses a flexible variant of the iterative adaptive-rank method developed in previous work by Ballani and Grasedyck \cite{ballani2013projection}. We combined this variant with flexible low-rank GMG preconditioners to enable efficient solution of both linear and nonlinear problems in high dimensions.

A central component of this work is the integration of adaptive-rank tensor methods with multigrid techniques in a matrix-free setting. By exploiting the separable structure of the underlying operators, the GMG hierarchy preserves low-rank structure across levels and provides an effective preconditioner for iterative methods. For nonlinear problems, this framework is embedded within an inexact Newton method, where the use of flexible GMG preconditioning enables robust and efficient solution of the resulting linear systems without requiring explicit Jacobian assembly.

The numerical experiments demonstrate the effectiveness of the proposed approach across several representative model problems. For the Poisson equation, the methods exhibit near mesh-independent convergence and scale favorably with dimension, while maintaining high compression rates. In the Dougherty–Fokker–Planck problem, the method remains robust in increasingly stiff regimes, with iteration counts that are largely independent of the discretization size and with low-rank structure that adapts naturally during relaxation toward equilibrium. For the Allen–Cahn equation, the method successfully captures complex interface dynamics while maintaining moderate ranks, and the inexact Newton solver demonstrates reliable convergence even during periods of rapid topological change. These results further highlight the capabilities of multilevel methods in creating robust preconditioners for high-dimensional linear and nonlinear PDEs.

Several promising directions remain for future work. From an algorithmic perspective, ongoing efforts include extending these methods to problems with variable coefficients and developing more sophisticated smoothers and coarse solvers to further improve robustness and efficiency. Although not considered here, adaptive time-stepping strategies may also be beneficial for strongly nonlinear or stiff problems. Another important direction is the application of these methods to asymptotic-preserving discretizations, which often involve implicit components and therefore require efficient solvers for high-dimensional systems. We are also exploring extensions to finite-element discretizations in order to treat more general geometries and boundary conditions. Finally, implementation on modern hardware is an important practical direction, with the potential to provide significant performance gains for large-scale high-dimensional simulations.

\section*{Acknowledgements}

This work was performed under the auspices of the U.S. Department of Energy by Lawrence Livermore National Laboratory under Contract DE-AC52-07NA27344. LLNL-POST-2015756; LLNL-JRNL-2020457. Sandia National Laboratories is a multimission laboratory managed and operated by National Technology \& Engineering Solutions of Sandia, LLC, a wholly owned subsidiary of Honeywell International Inc., for the U.S. Department of Energy’s National Nuclear Security Administration under contract DE-NA0003525. This paper describes objective technical results and analysis. Any subjective views or opinions that might be expressed in the paper do not necessarily represent the views of the U.S. Department of Energy or the United States Government. This work was partially supported by Sandia LDRD project \#233973. W.A.S. also acknowledges support from the U.S. Department of Energy through the Nicholas C. Metropolis Postdoctoral Fellowship under the Laboratory Directed Research and Development program at Los Alamos National Laboratory. Los Alamos National Laboratory report number LA-UR-26-24595.

\bibliographystyle{siam} 
\bibliography{ref}

\end{document}